\documentclass[11pt]{article}
\usepackage{graphicx, color}
\usepackage{amsmath,latexsym,amsbsy,amssymb}
\def\ds{\displaystyle}
\def\forall{\hbox{for all}~}
\def\L{{\bf L}}
\def\H{{\cal H}}
\def\D{{\cal D}}

\def\bft{{\bf t}}

\def\bfn{{\bf n}}
\def\bfe{{\bf e}}

\def\Supp{\hbox{Supp}}
\def\ve{\varepsilon}

\def\avint{-\!\!\!\!\!\!\int}

\def\G{{\cal G}}
\def\E{{\cal E}}
\def\I{{\cal I}}
\def\J{{\cal J}}

\def\R{I\!\!R}

\def\implies{\Longrightarrow}
\def\vp{\varphi}

\def\vs{\vskip 2em}

\def\v{\vskip 1em}
\def\O{{\cal O}}

\def\begi{\begin{itemize}}
\def\endi{\end{itemize}}
\def\C{{\cal C}}

\def\M{{\cal M}}
\def\T{{\cal T}}
\def\ov{\overline}
\def\Tilde{\widetilde}
\def\Hat{\widehat}
\def\bega{\begin{array}}
\def\enda{\end{array}}
\def\meas{\hbox{meas}}

\def\bel{\begin{equation}\label}
\def\eeq{\end{equation}}
\def\sqr#1#2{\vbox{\hrule height .#2pt
\hbox{\vrule width .#2pt height #1pt \kern #1pt
\vrule width .#2pt}\hrule height .#2pt }}
\def\square{\sqr74}
\def\endproof{\hphantom{MM}\hfill\llap{$\square$}\goodbreak}

\newtheorem{theorem}{Theorem}[section]
\newtheorem{lemma}[theorem]{Lemma}

\newtheorem{definition}[theorem]{Definition}
\newtheorem{remark}[theorem]{Remark}

\begin{document}
\title{\bf  Optimal Shapes for Tree Roots}
\vs
\author{Alberto Bressan$^{(*)}$ Sondre T.~Galtung$^{(**)}$,  and Qing Sun$^{(*)}$ \\~~\\
(*)~Department of Mathematics, Penn State University, \\
University Park,
Pa.~16802, USA.\\~~\\
(**)~Department of Mathematical Sciences, \\ NTNU -- Norwegian University of Science and Technology, \\ NO-7491 Trondheim, Norway.\\~~\\
E-mails: bressan@math.psu.edu,~sondre.galtung@ntnu.no,~qxs15@psu.edu.}
\maketitle

\begin{abstract} The paper studies a class  of variational problems, modeling
optimal shapes for tree roots.
Given a measure $\mu$ describing the distribution of root hair cells,
we seek to maximize a harvest functional $\H$, computing the total amount of water and nutrients
gathered by the roots, subject to a 
cost for transporting these nutrients from the roots to the trunk.
Earlier papers had established the existence of an optimal measure, and a priori bounds. Here we derive necessary conditions for optimality.   
Moreover, in space dimension $d=2$, we prove that the support of an optimal measure
is nowhere dense.
\end{abstract} 

\vs

\section{Introduction}
\label{s:1}
\setcounter{equation}{0}
Variational problems related to the optimal shape of tree roots were recently
considered in \cite{BPS, BSun}. 
Here one seeks an optimal measure $\mu$, describing the distribution of root hair cells.    The goal is to maximize a payoff, measuring the amount of 
water and nutrient absorbed by the roots, minus a cost for transporting these nutrients
to the base of the trunk.  
As in \cite{BCS},  given an open set $\Omega\subset\R^d$, the density of nutrients is modeled by the solution to 
an elliptic equation with measure coefficients. 
$$\Delta u + f(u) - u\mu~=~0,\qquad\qquad x\in\Omega,$$
with Neumann boundary conditions.
Here
$$\int u\, d\mu$$
yields the total harvest.
In addition, a ramified transportation cost is present. For a given $0<\alpha<1$, this is
described by the $\alpha$-irrigation cost 
of the measure $\mu$  from the origin~\cite{BCM, MMS, X03}.
The existence of an optimal measure was first proved in \cite{BSun}
under a constraint on the total mass of the measure $\mu$, and then in \cite{BPS}
in a more general setting.  

In the present paper we initiate a  study of the properties of these optimal measures.
Our first result provides 
necessary conditions for optimality. These take the form
\bel{nc0}\Phi(x)~=~c Z(x)\qquad\qquad\forall x\in \Supp(\mu),\eeq
for a suitable constant $c>0$.
Here $\Phi(x)$ measures the rate of increase of the total harvest, if the measure $\mu$
is locally increased at the point $x\in\Omega$. On the other hand, $Z(x)$ is the landscape function \cite{BCM,   S1}.  This is proportional to the rate of increase of the irrigation cost,
if the measure $\mu$ is locally increased at the point $x$.

In the second part of the paper we perform a detailed study of the equation (\ref{nc0}),
in dimension $d=2$.  Our main result, Theorem~\ref{t:41}, shows that 
the support of an optimal measure $\mu$ is nowhere dense.
To appreciate the physical meaning of this fact, one may observe that
water and nutrients can be moved around either by diffusion, or by ramified transport.
Diffusion comes for free,  but it is only effective at short distances.
A network of roots is thus needed to transport water and nutrients over 
longer distances,  while at small scales one can rely on diffusion alone.

The proof of Theorem~\ref{t:41} 
exploits the fact that the two functions $\Phi$ and $Z$ have very different
regularity properties, hence the set where they coincide must be small.
To help the reader, we outline here the main ideas.

\begin{figure}[ht]
\centerline{\hbox{\includegraphics[width=7cm]{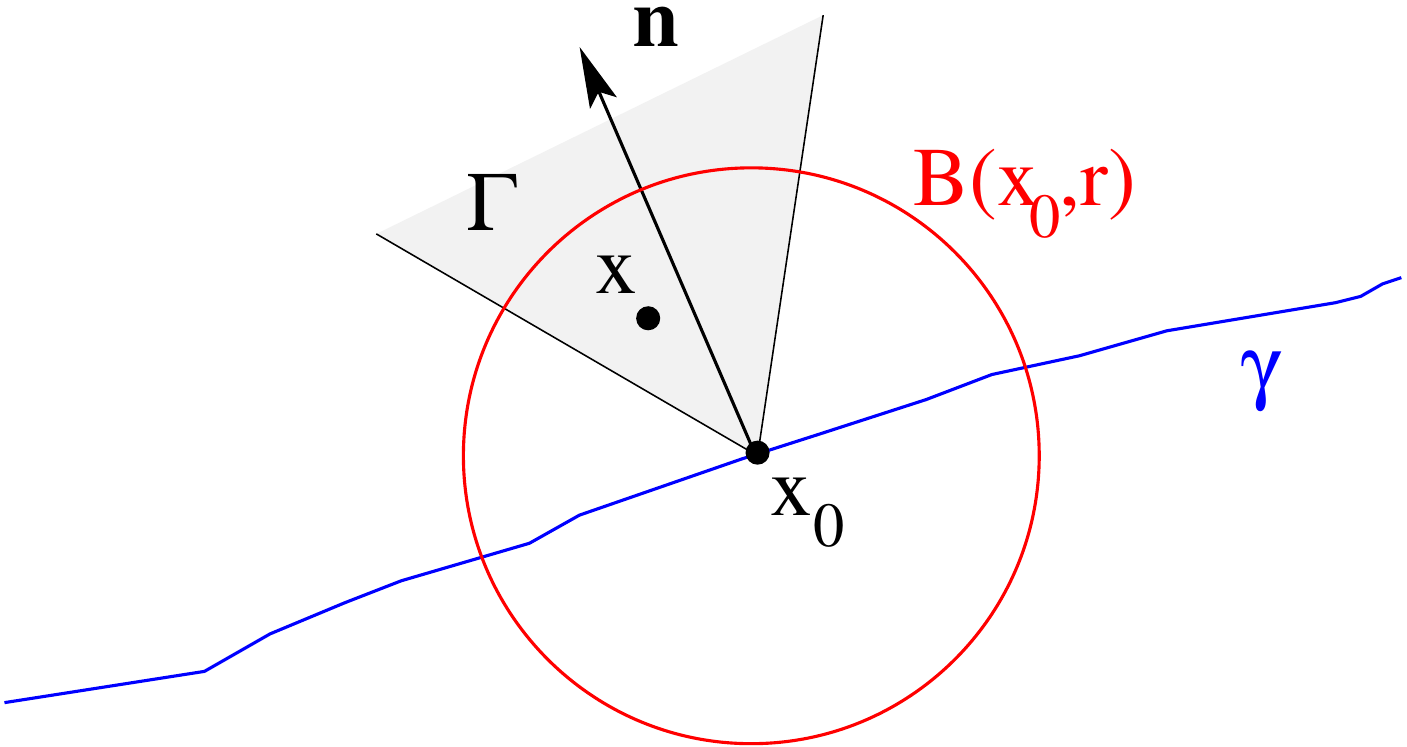}}}
\caption{\small Proving that $\Phi(x)<c Z(x)$, at several points $x$ near $x_0$.  }\label{f:ir156}
\end{figure}
Consider any particle path $s\mapsto \gamma(s)$, $s\in [0,\ov s]$, 
in an optimal irrigation plan for the measure $\mu$.     Given a point $x_0= \gamma(s_0)$,
let $\bfn$ be the unit  vector perpendicular to $\gamma$ at $x_0$, and consider the shaded region
in Fig.~\ref{f:ir156}
$$\Gamma~\doteq~\left\{ x\in\R^2\,;~\Big\langle \bfn\,,\, {x-x_0\over |x-x_0|}\Big\rangle~>~{1\over 2}\right\}.$$
By an argument based on Riesz' sunrise lemma~\cite{KF} we show that, for a.e.~$s_0\in [0,\ov s]$,
the landscape function $Z$ satisfies a lower H\"older estimate
on $\Gamma$, namely
\bel{01}Z(x)-Z(x_0)~\geq~\delta_0\cdot |x-x_0|^\alpha\qquad\quad\forall x\in \Gamma,\eeq
for some constant $\delta_0>0$ depending on $x_0$.

On the other hand, the function $\Phi$ can be bounded above by an
auxiliary function  $\Phi^+$, which satisfies
\bel{02}\left\{\bega{rll}\Delta \Phi^+&=~f\qquad &\hbox{on}~~\Omega,\\[3mm]
\Phi^+&\leq~c Z\qquad &\hbox{on}~~\gamma,
\enda\right.\eeq
for some source function $f$.
Since we only know that $f\in \L^1(\Omega)$,
from (\ref{02}) we do not obtain any 
pointwise upper bound on $\Phi^+$.
However, one can look at the 
average value of $\Phi^+$ over balls centered at $x_0$ with small radius $r>0$.
Relying on Vitali's covering theorem, together with an estimate on the 
Green's function
for the Laplacian on a suitable domain, we eventually obtain the
averaged integral estimate
\bel{03}
\avint_{\Gamma\cap B(x_0,r)} \Phi^+(x)\, dx- \Phi^+(x_0)~\leq~C_0\, r^{1-\ve},\eeq
where $\ve>0$ can be chosen arbitrarily small.  

Choosing $\ve< 1-\alpha$ and letting $r\to 0$, from (\ref{01})--(\ref{03}) we conclude that 
there exists a sequence of points $x_n\to x_0$ such that
$$\Phi(x_n)~\leq~\Phi^+(x_n)~<~c Z(x_n)\qquad\qquad\forall n\geq 1.$$
Since $\Phi$ is upper semicontinuous while $Z$ is lower semicontinuous,
this implies that the strict inequality $\Phi(x)< c Z(x)$ holds on an open set 
containing all the points $x_n$.  
Observing that the same conclusion can be reached for almost every point $x_0$ along
every irrigation path $\gamma$, this achieves the proof.

The remainder of the paper is organized as follows.   Section~\ref{s:2} reviews
the definition of the harvest functional~\cite{BCS, BPS, BSun}, and 
some  basic properties of
ramified transport and the landscape 
function~\cite{BCM, MMS, S1, X03}.
In Section~\ref{s:3} we formulate the optimization problem for tree roots, and 
derive a set of necessary  conditions for optimality, stated in 
Theorems~\ref{t:21} and \ref{t:22}. The proofs are an adaptation of the
arguments in~\cite{BSta}, where similar  necessary conditions were established for
a fishery model.    The key ideas are taken from~\cite{BDM}, where some
shape optimization problems involving measures were first studied.

The main new result of the paper, on the support of the optimal measure $\mu$, is stated in 
Theorem~\ref{t:41}.   The proof is worked out in Section~\ref{s:4}, 
relying on two key lemmas.   
Lemma~\ref{l:42}, providing a lower H\"older estimate on the landscape function, is then
proved in Section~\ref{s:5}. Finally, Lemma~\ref{l:43}, establishing an upper bound
on suitable averages of the function $\Phi^+$, is proved in Section~\ref{s:6}.

An alternative approach to the H\"older continuity of the landscape function
can be found in~\cite{BraSo11, BraSo14}. 
A different variational problem involving a ramified transportation cost has 
been recently studied  in~\cite{DS, PSX}.  Additional properties of optimal 
transportation plans were studied in~\cite{BCM1, BW, BS, DS07, MS, MoS}.
For a survey, see also~\cite{X15}.

\section{Review of the basic functionals}
\label{s:2}
\setcounter{equation}{0}

\subsection{Harvest functionals}
We consider a utility functional 
associated with plant roots.  Here the main goal is to collect
moisture and nutrients from the ground. 
To model the efficiency of a root, we let $u(x)$ be the density 
of water+nutrients at the point $x$, and consider a positive Radon measure $\mu$ 
describing the distribution  of root hair cells. 

To fix ideas, let  $\Omega\subset\R^d$, be an open set
of dimension $d\geq 2$, with $\C^2$ boundary.  We assume that 
$\mu$ is a positive, bounded Radon measure, supported 
 on the closure $\ov\Omega$, and absolutely continuous w.r.t.~capacity.
 This property is expressed by the implication
 \bel{cap0}\hbox{cap}_2(V)~=~0\qquad\implies\qquad \mu(V)~=~0.\eeq
 For the definition and basic properties of capacity we refer to \cite{AG, EG}.
Based on physical considerations, following \cite{BCS, BPS, BSun} we
consider the solution to   the elliptic problem with measure source
\bel{emu} \Delta u + f(u) - u\,\mu~=~0\qquad\qquad\hbox{on}~\Omega,\eeq
and  Neumann boundary conditions
\bel{NBC}
 \partial_{\bfn(x)}  u~=~0\qquad \hbox{on}~~\partial\Omega\,.\eeq
Here $\bfn(x)$ denotes the unit outer normal vector
at the boundary point $x\in\partial\Omega$, while
$\partial_\bfn u$ is the derivative of $u$ in the
normal direction.  Suitable assumptions on the source function $f$ will be given at (\ref{fp1}).

Elliptic problems with measure data have been 
studied in various papers \cite{BG1, BGO, DMOP},
and are now well understood. 
A key fact is that, roughly speaking, the Laplace operator ``does not see" sets with zero capacity. Therefore, a measure concentrated on a set with zero capacity 
does not affect the solution of (\ref{emu}).
Following \cite{BG1, BGO},  we denote by $\M_0$ the family of all bounded Radon 
measures which vanish on sets with zero capacity.

\begin{definition} Let $\mu$ be a positive, bounded Radon measure on $\ov\Omega$, 
which is absolutely continuous w.r.t.~capacity.
 A function
$ u\in\L^\infty(\Omega)\cap H^1(\Omega)$ is a solution to  
the elliptic problem (\ref{emu})-(\ref{NBC}) if
 \bel{nsol}
 -\int_\Omega \nabla u\cdot \nabla\vp\, dx+\int_\Omega f(u)\vp\, dx -\int_{\ov\Omega} u\vp\, d\mu
 ~=~0\eeq
 for every test function $\vp\in \C^\infty_c(\R^d)$.
\end{definition}

\begin{definition} 
In connection with a solution $u$ of (\ref{emu})-(\ref{NBC}),
the {\bf total harvest}  is defined as
\bel{HF1}\H(u,\mu)~\doteq~\int_{\ov\Omega}~ u\, d\mu\,.\eeq
\end{definition}
For reader's convenience we collect the main assumptions used throughout the paper.
\begi
\item[{\bf (A1)}]  {\it ~$\Omega\subset\R^d$ is a  bounded, connected open set with $\C^2$ boundary. Moreover, $0\notin \ov\Omega$.}
\item[{\bf (A2)}] {\it $f:\R\mapsto\R$ is a $\C^2$ function such that, for 
some constants $u_{max},K>0$,
\bel{fp1} f(u_{max})=0,\qquad \quad 
0\leq
f(u)\leq K,\qquad f''(u)< 0 \quad\forall u\in [0,u_{max}].\eeq}
\item[{\bf (A3)}] {\it The space dimension is $d\geq 2$. The exponent $\alpha$ in the 
irrigation cost satisfies \bel{alp}\alpha~>~ 1-{1\over d-1}\,.\eeq
}
\endi

\begin{remark}{\rm
   If $\mu$ is a general measure
and $u$ is a discontinuous function,  the integral (\ref{HF1}) may not be well defined.
To resolve this issue,  calling
 $$\avint_V u\, dx~=~{1\over \hbox{meas}(V)}\,\int_V u\, dx$$
the average value of $u$ on a set $V$, for each $x\in \ov\Omega$ we consider the limit
\bel{upv}
u(x)~=~\lim_{r\downarrow 0}~\avint_{\Omega\cap B(x,r)} u(y)\, dy.\eeq
As proved in  \cite{FZ}, if $u\in H^1(\Omega)$ then
the above limit exists at all points $x\in \ov\Omega$ with the possible exception of a set 
whose capacity is zero. If  $\mu\in \M_0$, then the integral (\ref{HF1}) is  well defined.
Our present setting is actually even better, because in (\ref{emu}) $u$ and $\mu$ are positive while $f$ is bounded.
Therefore, if the constant $C$ is chosen large enough, the function $u+C|x|^2$ is subharmonic \cite{AG}.  As a consequence, 
the limit (\ref{upv}) is well defined at every point $x\in \ov\Omega$.
}\end{remark}

\v
\subsection{Optimal irrigation plans}

Given $\alpha\in [0,1]$ and a positive measure $\mu$ on $\R^d$, the minimum
cost for irrigating the measure $\mu$ from the  
origin will be denoted by $\I^\alpha(\mu)$.
Following 
\cite{MMS}, this cost
 can be defined as follows.
Let $M=\mu(\R^d)$ be the total mass to be transported and let 
$\Theta=[0,M]$.
We think of each $\theta\in \Theta$ as a ``water particle".

\begin{definition}
A measurable map 
\bel{iplan}
\chi:\Theta\times \R_+~\mapsto~ \R^d\eeq
is called an {\bf admissible irrigation plan} for the measure $\mu$
if 
\begi
\item[(i)] For every $\theta\in \Theta$, the map
$t\mapsto \chi(\theta,t)$ is 1-Lipschitz. 
More precisely, for each $\theta$ there exists a stopping time $T(\theta)$ such that, calling 
$$\dot \chi(\theta,t)~=~{\partial\over\partial t} ~\chi(\theta,t)$$
the partial derivative w.r.t.~time, one has
\bel{stime}\bigl|\dot \chi(\theta,t)\bigr|~=~\left\{ \bega{rl} 1\qquad &\quad\hbox{for a.e.}
~t\in [0, T(\theta)],\\[3mm]
0\qquad &\quad\hbox{for}
~t\geq T(\theta).\enda\right.\eeq
\item[(ii)] At time $t=0$ all particles are at the origin:
$$\chi(\theta,0)~=~{\bf 0}~\in~\R^d \qquad\forall \quad \theta\in\Theta.$$
\item[(iii)] The push-forward of the Lebesgue measure on $\Theta=[0,M]$ through the map $\theta\mapsto 
\chi(\theta,T(\theta))$ coincides with the measure $\mu$.
In other words, for every open set $A\subset\R^n$ there holds
\bel{chi1}\mu(A)~=~\hbox{meas}\Big( 
\bigl\{ \theta\in \Theta\,;~~\chi(\theta,T(\theta))\in A\bigr\}\Big),\eeq
where $\meas(\cdot)$ denotes the Lebesgue measure on $\Theta$.
\endi
\end{definition}

In order to define the corresponding transportation cost, we first consider the amount of 
paths which go through a point $x\in \R^d$:
\bel{chi}|x|_\chi~=~\meas\Big(\bigl\{\theta\in \Theta\,;~~\chi(\theta,t)= x~~~\hbox{for some}~~t\geq 0\bigr\}\Big).\eeq
We think of $|x|_\chi$ as the {\bf total flux through the
point $x$}.
\v
\begin{definition} 
For a given $\alpha\in [0,1]$,
the total cost of the irrigation plan $\chi$ is
\bel{TCg}
\E^\alpha(\chi)~\doteq~\int_\Theta\left(\int_{\R_+} \bigl|\chi(\theta,t)
\bigr|_\chi^{\alpha-1} \cdot |\dot \chi(\theta,t)|\, dt\right)
d\theta.\eeq
The  {\bf $\alpha$-irrigation cost} of a measure $\mu$
is defined as 
\bel{Idef}\I^\alpha(\mu)~\doteq~\inf_\chi \E^\alpha(\chi),\eeq
where the infimum is taken over all admissible irrigation plans.

We say that $\mu$ is $\alpha$-irrigable if $\I^\alpha(\mu)<+\infty$.
\end{definition}

A lower bound on the transportation cost is provided by

\begin{lemma}
For any positive Radon measure $\mu$ on  $\R^d$ and any $\alpha\in [0,1]$, one has
 \bel{Ilb}
\I^\alpha(\mu)~\geq~\int_0^{+\infty}\Big[\mu\bigl(\{x\in\R^d;~|x|\geq r\}
\bigr)\Big]^\alpha dr\,.\eeq
In particular, for every $r>0$ one has
 \bel{Irb}
\I^\alpha(\mu)~\geq~r\cdot \Big[\mu\bigl(\{x\in\R^d;~|x|\geq r\}
\bigr)\Big]^\alpha.\eeq
\end{lemma}

We recall that optimal irrigation plans
satisfy
\v
{\bf Single Path Property:} {\it  If $\chi(\theta, \tau)=\chi(\theta',\tau')$ for some 
$\theta, \theta'\in\Theta$ and 
$0<\tau\leq \tau'$, then 
\bel{SPP}\chi(\theta, t)~=~\chi(\theta', t)\qquad\forall t\in [0, \tau].\eeq
}

\v
\begin{remark} {\rm In the case $\alpha=1$, the expression (\ref{TCg}) reduces to
$$
\E^\alpha(\chi)~\doteq~\int_\Theta\left(\int_{\R_+} |\dot \chi_t(\theta,t)|\, dt\right)
d\theta~=~\int_\Theta[\hbox{total length of the path} ~\chi(\theta,\cdot)]\, d\theta\,.$$
Of course, this length is minimal if every path $\chi(\cdot,\theta)$
is a straight line, joining the origin with $\chi(\theta, T(\theta))$.  Hence
$$\I^\alpha(\mu)~\doteq~\inf_\chi \E^\alpha(\chi)~=~\int_\Theta |\chi(\theta, T(\theta))|\, d\theta~=~\int |x|\, d\mu\,.$$

On the other hand, when $\alpha<1$, moving along a path which is traveled by few other particles
comes at a  high cost. Indeed, in this case the factor $\bigl|\chi(\theta,t)
\bigr|_\chi^{\alpha-1}$ becomes  large.   To reduce the total cost,  is thus convenient
that particles travel along the same path as far as possible.}
\end{remark}

%
\v
\subsection{The landscape function}
Let $\chi:\Theta\times\R_+\mapsto\R^d$ be a
(possibly not optimal) irrigation plan which satisfies the single path property.
The  {\bf landscape function} $Z:\R^d\mapsto\R_+$ is defined as follows \cite{BCM, S1}.

\begi
\item[(i)] As a first step, consider the  set 
\bel{Gamdef}\Gamma~\doteq~\{x\in\R^d\,;~|x|_\chi>0\}.\eeq 
We think of $\Gamma$ as the union of all irrigation paths.
For $x\in\Gamma$, choose any particle $\theta\in \Theta$ that reaches $x$, so that
$\chi(\theta,\tau)=x$ at some time $\tau\geq 0$.  We then define
\bel{Zdef}
Z_\chi(x)~\doteq~\int_0^\tau \bigl|\chi(\theta,t)\bigr|_\chi^{\alpha-1}\, dt.\eeq
By the single path property (\ref{SPP}), the above integral is independent of the choice of $\theta$.
\item[(ii)] The landscape function $Z$ is now defined as 
the lower semicontinuous envelope of $Z_\chi$. Namely,
\bel{Zd2}Z(x)~\doteq~\liminf_{y\in \Gamma, y\to x} Z_\chi(y),\eeq
with the understanding that  $Z(x)=+\infty$ if $x\notin \ov\Gamma$.
\endi
	
The following results were originally proved for a probability measure. 
However, by a rescaling, it is clear that they remain valid for any bounded positive measure. 
%
	\begin{lemma}
Let $\mu$ be a $\alpha$-irrigable measure on $\R^d$, and 
let $\chi$ be an optimal irrigation plan for $\mu$. 
Calling $Z$ the landscape function of $\chi$, one has
	\bel{4x0} \I^\alpha(\mu)~=~\int_{\R^d} Z(x)\, d\mu\, . \eeq
\end{lemma}
{\bf Proof.}  See Corollary 4.5 in \cite{S1}.
%
\begin{theorem}
	Let $\mu$ be an $\alpha$-irrigable measure on $\R^d$ and 
	let $g$ be a measurable function such that $\|g\|_{L^\infty(\mu)}\leq 1$. 
	Consider the measure $\nu \doteq (1 + g)\mu$. Then
	\bel{4x1} \I^\alpha(\nu)~\leq~\I^\alpha(\mu) + \alpha \int_{\R^d} Z(x)g(x)\, d\mu\, .  \eeq	
\end{theorem}
{\bf Proof. } See Theorem 4.7 in \cite{S1}. Notice that the condition $\|g\|_{L^\infty(\mu)}\leq 1$ guarantees that $\nu$ is a positive measure. 

\begin{lemma}\label{l:LH}
	Let $\chi$ be an optimal $\alpha$-irrigation plan for the measure $\mu$, and let 
	$Z(\cdot)$ be the corresponding landscape function. Then, for any two point $x\in \Gamma, y\in \Gamma$,  one has
	\bel{0xa0}
	Z(x) - Z(y) ~\leq~{1\over \alpha} |x|_\chi^{\alpha-1}|x - y|\, .
\eeq	
\end{lemma}
{\bf Proof. } See Corollary 3.10 in \cite{BraSo11}.
%
%
%

We conclude this section by observing that, for any  branch in $\Gamma$, the 
arc-length can be bounded below in terms of the Euclidean distance. 

\begin{lemma}\label{l:112} Let $\chi$ be an optimal irrigation plan  for the measure $\mu$.
Consider any path $t\mapsto \gamma(t) = \chi(\theta,t)$ for some particle
$\theta\in \Theta$.   Assume that its multiplicity is bounded from below:
\bel{mbb}
|\gamma(s)|_\chi~\geq~\delta_0~>~0\qquad\qquad\forall s\in [0,\ov s].\eeq
Then there exists a constant $C$ such that 
\bel{distlb} |t-s|~\leq~
C\,|\gamma(t)-\gamma(s)|\qquad\qquad\forall s,t\in [0,\ov s].\eeq
\end{lemma}
\v
{\bf Proof.}  To fix ideas, assume that  $s<t$,  $\,x=\gamma(s)$, $\,y=\gamma(t)$.
By optimality,  the multiplicity function $\tau\mapsto |\gamma(\tau)|_\chi$ is non-increasing along $\gamma$. Therefore
\bel{0xd5} Z(y) - Z(x)~=~\int_s^t |\gamma(\tau)|_\chi^{\alpha-1}\, d\tau
~\geq~|\gamma(s)|_\chi^{\alpha-1}(t-s)\, . \eeq

On the other hand, according to Lemma \ref{l:LH} we have
\bel{0xd6} Z(y) - Z(x)~\leq~{1\over \alpha} |\gamma(t)|_\chi^{\alpha-1}|x - y|\, . \eeq
Calling $M$ the total mass of the irrigated measure,
combining (\ref{0xd5}) with (\ref{0xd6}) one obtains
\bel{0xd7} t-s~\leq~{1\over\alpha} \left({|\gamma(t)|_\chi
\over |\gamma(s)|_\chi}\right)^{\alpha-1}\, |y-x|~\leq~{1\over\alpha} \left({\delta_0\over M}\right)^{\alpha-1}\, \bigl|\gamma(t)-\gamma(s)\bigr|.\eeq
\endproof

\v
\section{Necessary conditions for optimal tree roots}
\label{s:3}
\setcounter{equation}{0}

Following \cite{BPS, BSun}, the optimization problem for tree roots can be stated as
\begi
\item[{\bf(OPR)}]  {\it Maximize the functional 
\bel{mH} \H(u,\mu) - c
\I^\alpha(\mu),\eeq
among all couples $(u,\mu)$, where $\mu$ is a positive measure on $\ov\Omega$, and $u$ is a solution to (\ref{emu})-(\ref{NBC}).}
\endi
Existence of solutions was proved in \cite{BPS}. 
\begin{theorem}\label{t:20} Let the assumptions {\bf (A1)--(A3)} hold. Then the problem 
{\bf (OPR)} has at least one optimal solution $(u^*,\mu)$, 
satisfying
\bel{harveq}\left\{ \bega{rll}
\Delta u^* + f(u^*) -  u^* \mu &=~0\qquad\qquad &x\in\Omega\,,\\[2mm]
\bfn \cdot \nabla u^* &=~0\qquad\qquad &x\in\partial\Omega\,.\enda\right.\eeq
The  measure $\mu$ on 
$\ov\Omega$ has bounded total mass.
\end{theorem}

Indeed, the result in \cite{BPS} established the existence of an optimal pair $(u^*,\mu)$,
in the more general case where the set $\Omega\subseteq\R^d$ may be unbounded,
possibly also with $0\in \ov\Omega$.  
The analysis in \cite{BPS} also shows  that the irrigation cost
for the optimal measure is bounded: $\I^\alpha(\mu)< +\infty$.
By itself, this does not guarantee that the total mass of the
measure $\mu$ is bounded, because $\mu$ may 
concentrate an infinite amount of mass near the origin, 
where the transportation cost is almost zero.
In the present setting however, 
thanks to the additional assumption $0\notin \ov\Omega$ 
in {\bf (A1)},  by (\ref{Irb}) we conclude that the total mass of $\mu$ is bounded by
$$\mu(\ov\Omega)~\leq~\left( \I^\alpha(\mu)\over r_0\right)^{1/\alpha},\qquad\qquad 
r_0\,\doteq\, d(0, \ov \Omega) \,=\, \min\bigl\{ |x|\,;~x\in \ov\Omega\bigr\}.$$
\v

The main  goal of this section
is to derive necessary conditions for optimality. 

\begin{theorem}\label{t:21} Let the assumptions {\bf (A1)--(A3)} hold.
Let  $(u^*,\mu)$ be an optimal solution to the problem {\bf (OPR)}, satisfying 
(\ref{harveq}).
Let $\chi$ be an optimal irrigation plan for the measure $\mu$, and 
let $Z$ be the corresponding landscape function.

Then there exists a  bounded solution $\psi\geq 0$ to the adjoint equation
\bel{adj}\left\{ \bega{rll}
\Delta \psi + f'(u^*)\psi - \psi \mu
&= ~-\mu\qquad\qquad &x\in\Omega\,,\\[2mm]
\bfn \cdot \nabla\psi &=~0\qquad\qquad &x\in\partial\Omega\,,\enda
\right.\eeq
such that,  $\mu$-almost everywhere, one has
\bel{NC3}
(1-\psi)u^*~ =~c\,\alpha Z.\eeq
\end{theorem}
\v
{\bf Proof.}  We follow the same steps as in the proof of Theorem~2.1 in \cite{BSta},
with minor modifications.
\v
{\bf 1.} We begin by proving that the solution $u^*$ of (\ref{harveq}) is uniformly positive on 
$\ov\Omega$.  Indeed, since $0\notin \ov\Omega$, recalling (\ref{0xd5}) we can choose 
a constant $0<\delta_0< u_{max}$ such that
the landscape function satisfies
\bel{Z>}
c\alpha\,Z(x)~\geq~\delta_0~>~0\qquad\qquad\forall x\in \ov\Omega.\eeq
We now claim that the optimal measure $\mu$ vanishes on the set where $u^* < c\alpha Z$,
namely
\bel{mu00}
\mu\Big(\bigl\{ x\in \ov\Omega\,;~~u^*(x)\,<\, c\alpha Z(x)\bigr\}\Big)~=~0.
\eeq
Indeed, if (\ref{mu00}) fails, we can consider the reduced measure 
$\mu_0\,\doteq\,g\,\mu$, where
$$g(x)~=~\left\{\bega{rl} 1\quad &\hbox{if}~~u^*(x)\geq c\alpha Z(x),\\[1mm]
0\quad &\hbox{if}~~u^*(x)< c\alpha Z(x).\enda\right.$$
Let $u_0$ be the solution to 
\bel{ha0}\left\{ \bega{rll}
\Delta u + f(u) -  u \mu_0 &=~0\qquad\qquad &x\in\Omega\,,\\[2mm]
\bfn \cdot \nabla u &=~0\qquad\qquad &x\in\partial\Omega\,.\enda\right.\eeq
In view of the assumptions (\ref{fp1}) on the source function $f$, the existence of 
such a solution follows from the analysis in \cite{BCS}.   Moreover, a comparison
argument yields
$$0~\leq~u_0(x)~\leq~u_{max}\qquad\qquad\forall x\in\Omega.$$
Since $u^*$ provides a subsolution to (\ref{ha0}),
we have $u^*\leq u_0$.   Hence, by (\ref{4x1}),
$$\bega{l} \H(u_0, \mu_0) -  \H(u^*, \mu)- c\Big[\I^\alpha(\mu_0)- \I^\alpha(\mu)
\Big] ~ =~\ds \int_{\ov\Omega}u_0\,d\mu_0  - \int_{\ov\Omega}u^*\,d\mu - c\Big[\I^\alpha(\mu_0)- \I^\alpha(\mu)
\Big] \\[4mm]
\qquad\qquad \ds \geq~ \int_{\ov\Omega} (g-1) u^*\,d\mu 
-c \alpha \int_{\ov\Omega} (g-1) Z \, d\mu~>~0,\enda $$
against  the optimality
of $(u^*,\mu)$.  Therefore,  $\mu_0=\mu$ and (\ref{mu00}) holds.

Next,
in view of (\ref{mu00}), the function
$$\Tilde u(x)~\doteq~\max\bigl\{ \delta_0\, ,~u^*(x)\bigr\}$$
is a subsolution of (\ref{ha0}). Indeed, on the set where $\Tilde u(x) = \delta_0$
we have
$$\Delta \Tilde u + f(\Tilde u) - \Tilde u\,\mu~=~f(\Tilde u) ~=~f(\delta_0)~>~0.$$
On the other hand, by (\ref{fp1}) the constant function $u(x)=u_{max}$ is trivially a supersolution.
We thus conclude that 
\bel{Tus}
u_{max}~\geq~u^*(x) ~\geq ~\Tilde u(x)~\geq ~\delta_0~>~0\qquad\quad\forall x\in\ov\Omega\,.\eeq
\v
{\bf 2.}
Consider a family of perturbed measures, of the form
\bel{mep}\mu_\ve~=~\mu+\ve \nu,\eeq
where 
\bel{nudef}  \nu ~=~ g\,\mu,\qquad \hbox{with}\quad
\|g\|_{\L^\infty} \leq 1.\eeq
Let $u_\ve$  
be the corresponding solution of 
\bel{ueps}\Delta u_\ve +f(u_\ve) -
u_\ve\, \mu_\ve~=~0\qquad \qquad \hbox{on}~~ \ov\Omega,\eeq
with Neumann boundary conditions.

 When the measure $\mu$ is replaced by $\mu_\ve$,  
 by (\ref{4x1}) the irrigation cost satisfies
\bel{Iae}
\I^\alpha(\mu_\ve)~\leq~\I^\alpha(\mu) + \alpha \ve \int_{\ov\Omega} Z(x)g(x) \, d\mu(x).\eeq

In the next steps we shall derive a formula computing the corresponding change in the harvest functional $\H(u_\ve, \mu_\ve)$.

%
%
\v
{\bf 3.} Following \cite{BSta, BDM}, consider the space 
 $X_\mu\doteq \mathbf{H}^1(\Omega)\cap \mathbf{L}^2(\mu)$, and its dual space
 $X'_\mu$. As shown in \cite{BSta}, assuming that $\mu$ is not the zero measure,
the space $X_\mu$ is a Hilbert space with the equivalent inner product
\bel{ipmu}\langle u,v\rangle_{X_\mu}\doteq \int_{\Omega} Du\cdot Dv~dx+\int_{\ov\Omega} uv~d\mu\, .\eeq
Given $F\in X_\mu'$,  consider the problem of finding $u\in X_\mu$
which satisfies
\bel{lin}
 \Delta u~-~u\mu~=~F,\eeq
with Neumann boundary conditions (\ref{NBC}).  
We define the resolvent operator $R_\mu: X^{'}_\mu(\Omega)\to X_\mu(\Omega)$ 
by setting $R_\mu(F)=u$, where $u$ is the unique solution of (\ref{lin}). 
By Riesz' theorem, 
$R_\mu$ is a bounded linear operator from $X^{'}_\mu(\Omega)$ onto $X_\mu(\Omega)$, and thus continuously differentiable.
\v
{\bf 4.} 
Now let $\mu_{\ve}$, $u_\ve$  be as in (\ref{mep}), (\ref{ueps}).
%
%
Notice that (\ref{ueps})
is equivalent to
\bel{neff}
\Delta u_\ve~-~u_\ve\mu~=~-f(u_\ve)+\ve  u_\ve \nu .
\eeq
 Using the resolvent operator, (\ref{neff}) can be written as 
 \bel{ueR}u_\ve~=~R_{\mu}\bigl(-f(u_\ve)+\ve u_\ve  \nu\bigr)\, .\eeq
To prove that the map $\ve\mapsto u_\ve$ is differentiable,
consider the function $\Psi: \R\times X_{\mu}\to X_{\mu}$ defined as
\bel{impl}
\Psi(\ve,w)~\doteq~ w-R_{\mu}\Big(-f(w)+\ve w\, \nu\Big).
\eeq
Being the composition of the linear operator $R_{\mu}$ and a smooth map,
it is clear that $\Psi$ is continuously differentiable.  When $\ve=0$ we already know that
\bel{old}\Psi(0, u^*) ~=~0~\in ~X_\mu\,.\eeq
We claim that, for $\ve$ in a neighborhood of zero, the equation 
\bel{old2}\Psi(\ve, w) ~=~0\eeq
implicitly defines a function $w(\ve)=u_\ve$, providing the solution to (\ref{neff}). 

As shown in step {\bf 6} of the proof of Theorem~2.1
in \cite{BSta},
 the linear operator
\bel{op3}w~\mapsto~w+R_\mu\bigl(f'(u^*) w\bigr)\eeq
has a bounded inverse on $X_\mu$. By the implicit function theorem, it follows that 
the map $\ve\mapsto u_\ve$ is well defined, and differentiable in a 
neighborhood of the origin.

%
%
%
  Having established the differentiability of the map $\ve\mapsto u_\ve$, 
 its derivative at $\ve=0$ can be computed by differentiating (\ref{ueR}).  This yields
$$v~\doteq~\frac{d u_{\ve}}{d\ve}\Big|_{\ve=0}~=~
R_\mu\bigl(-f'(u^*) v +  u^*\nu \bigr).$$
Therefore $v$ satisfies the
linear, non-homogeneous equation
\bel{41}\left\{ \bega{rll}
\Delta v + f'(u^*) v -  v \mu &
= ~u^* \nu\qquad\qquad &x\in\Omega\,,\\[2mm]
\bfn \cdot \nabla v &=~0\qquad\qquad &x\in\partial\Omega\,.\enda\right.\eeq
Notice that (\ref{41}) could be formally obtained by inserting the expansion
$$u_\ve~=~u^*+\ve v + o(\ve)$$
in (\ref{neff}), and retaining terms of order $\O(\ve)$.
\v
{\bf 5.} In this step we show that the adjoint problem (\ref{adj}) has a uniformly 
bounded solution $\psi\in X_\mu$. 
Toward this goal, we first choose $\lambda>0$ large enough so that 
\bel{ff'}
f'(u) (\lambda u +1)~<~\lambda f(u)\qquad \forall  u\in [\delta_0,  u_{max}].\eeq
Notice that such a constant exist, thanks to (\ref{Tus}) and  the assumptions (\ref{fp1}).

We now
claim that, the function
\bel{psi+} \psi^+~=~\lambda u^* + 1\eeq
is a supersolution to (\ref{adj}).
Indeed, inserting (\ref{psi+}) in (\ref{adj}) and using (\ref{harveq}), we obtain
$$\bega{rl}
\Delta \psi^+ + f'(u^*)\psi^+ + (1-\psi^+) \mu &= ~\lambda \Delta u^* + f'(u^*) (\lambda u^*+1)
 - \lambda u^*\mu\\[2mm]
&=~f'(u^*)(\lambda u^* + 1) - \lambda f(u^*)~\leq~0.\enda
$$
This holds for every $x\in \ov\Omega$, because of (\ref{Tus}) and (\ref{ff'}).
We thus conclude that $\psi$ satisfies the uniform bounds
\bel{psib}
0~\leq~\psi(x)~\leq~\lambda \,u^*(x) + 1~\leq~\lambda \,u_{max} + 1.\eeq
\v
{\bf 6.} 
Next, let $\psi$ be the solution to the adjoint problem (\ref{adj}). Using $v$ as test function
and integrating by parts one obtains
\bel{int3}\int v\, d\mu~=~\int \nabla \psi\cdot \nabla v\, dx - \int f'(u^*)  \psi v\, dx + \int \psi v\, d\mu
~=~-\int \psi u^* d\nu.\eeq
Notice that the last identlity follows from the fact that $v$ is a weak solution to (\ref{41}),
using $\psi$ as test function.

Differentiating the harvest functional  w.r.t.~$\ve$ and using (\ref{int3})
one obtains
\bel{42}{d\over d\ve}\H(u_\ve, \mu_\ve)\bigg|_{\ve=0}~=~
\int_{\ov \Omega} u^*\,d\nu +
\int_{\ov\Omega} v\,d\mu~=~\int_{\ov\Omega} (1-\psi) u^*\, d\nu  .\eeq
\v
{\bf 7.} 
Since $(u^*, \mu)$ yield an optimal solution, in view of (\ref{42}),  (\ref{nudef}), and (\ref{Iae}), we obtain
\bel{necc}
\bega{rl}0&\geq~\ds
\limsup_{\ve\to 0+}\,  \left[{\H(u_\ve,\mu_{\ve})-\H(u^*,\mu)
\over\ve}-c\,{ \I^\alpha(\mu_{\ve})-\I^\alpha(\mu)\over\ve}\right]
\\[4mm]
&\geq~\ds
\int_{\ov\Omega}\Big((1-\psi)u^*-c \alpha\, Z\Big)g\,d\mu\,.\enda
\eeq
Since the function $g\in \L^\infty$ can be chosen arbitrarily, we conclude that 
the identity (\ref{NC3}) must hold almost everywhere w.r.t.~the measure $\mu$.
\endproof
\vs

Outside the support of $\mu$, 
the identity (\ref{NC3}) may fail.  Yet, we expect that it can be replaced by an inequality.
A result in this direction, valid in dimension $d=2$, is now proved.

\begin{theorem}\label{t:22} Assume $d=2$.  In the same setting as Theorem~\ref{t:21},
let  $(u^*,\mu)$ be an optimal solution to the problem {\bf (OPR)}, and let $\chi$ be an optimal irrigation
plan for the measure $\mu$. Consider any   particle path
$s\mapsto\gamma(s)\doteq \chi(\bar \theta, s)$, $s\in [0,\bar s]$, 
where the multiplicity remains strictly positive.
Then,
\bel{49}
(1-\psi)u^*~ \leq~c\,\alpha Z.\eeq
at almost every point $x=\gamma(s)$, $s\in [0, \ov s]$, such that $\gamma(s)\in \Omega$.
\end{theorem}

{\bf Proof.} {\bf 1.} Assume that the conclusion does not hold. Then the 
set 
$$S~\doteq~\Big\{s\in [0, \ov s]\,;~~\gamma(s)\in\Omega,~~\bigl(1-\psi(\gamma(s)\bigr)u^*\bigl(\gamma(s)\bigr)
~ >~c\,\alpha Z\bigl(\gamma(s)\bigr)\Big\},$$
where the inequality (\ref{49}) fails,
has positive Lebesgue measure.   

Let $\nu$ be the  measure supported along the 1-dimensional curve $\gamma$, 
obtained as the push-forward of Lebesgue measure on $S$, via the map
$s\mapsto \gamma(s)$. 

For $\ve>0$, consider the measures
$$\mu_\ve ~\doteq~ \mu+\ve \nu,$$ and let $u_\ve$ be the corresponding solutions to 
(\ref{emu})-(\ref{NBC}).
Since the dimension is $d=2$, in view of Lemma~\ref{l:112} it follows
that the measure $\nu$ is absolutely continuous w.r.t.~capacity.   
By a similar argument as in the proof of the previous theorem, 
the derivative of the harvest functional is computed by
\bel{dHe} {d\over d\ve} \H(u_\ve, \mu_\ve)\bigg|_{\ve=0}~=~\int_{\Omega} 
(1-\psi) u^*\, d\nu\,.\eeq
\v
{\bf 2.} It remains to estimate the change in the irrigation cost.
The measure $\mu_\ve$ has total mass 
$$\mu_\ve(\ov\Omega)~=~\mu(\ov\Omega) + \ve \nu(\Omega)~=~M+\ve\,\meas(S),$$
where $\meas(\cdot)$ always denotes Lebesgue measure.
Since the measure $\nu$ is supported along the curve $\gamma$, it is natural to 
consider an irrigation plan $$\chi_\ve: [0, M+\ve\,\meas(S)]\times \R_+~\mapsto~\R^2$$
which coincides with $\chi$ for $\theta\in [0,M]$,  while all the additional particles 
$\theta\in\, \bigl]M, M+\ve\, \meas(S)
\bigr]$
are transported to destination along the same path $\gamma$.
The change in multiplicity at points $x=\gamma(s)$ is thus 
$$\bigl| \gamma(s)\bigr|_{\chi_\ve}~=~\bigl| \gamma(s)\bigr|_\chi + \ve \,\meas\bigl(
S\cap [s, \ov s]\bigr).$$
In turn, the increase in the irrigation cost is computed as
\bel{icep}\bega{rl} \E^\alpha(\chi_\ve) - \E^\alpha(\chi)&\ds
=~\int_0^{\ov s} \Big( \bigl| \gamma(s)\bigr|_{\chi_\ve}^\alpha 
-\bigl| \gamma(s)\bigr|_{\chi}^\alpha\Big) \, ds\\[4mm]
&\ds =~\ve\, \alpha \int_0^{\ov s} \bigl| \gamma(s)\bigr|_{\chi}^{\alpha-1}
\meas(S\cap [s, \ov s]\bigr) ds+o(\ve)\\[4mm]
&\ds =~\ve\, \alpha \int_S Z(\gamma(s))\, ds + o(\ve).\enda\eeq
Here the last identity follows from 
$$ \bega{rl}\ds
\int_0^{\ov s} \bigl| \gamma(s)\bigr|_{\chi}^{\alpha-1}
\meas(S\cap [s, \ov s]\bigr) ds &=~\ds \int_0^{\ov s} \Bigl( {d \over ds} \int_0^s \bigl| \gamma(t)\bigr|_{\chi}^{\alpha-1} \,dt \Bigr)
\Bigl( \int_s^{\ov s} {\bf 1}_S(t)\,dt \Bigr)\,ds \\[4mm]
&= \ds\Bigl[ Z(s) \int_s^{\ov s} {\bf 1}_S(t)\,dt \Bigr]_0^{\ov s}-\int_0^{\ov s} Z(s) (-{\bf 1}_S(s))\,ds,
\enda
$$
where we have used the fact that $Z(0) = 0$.

Since $\chi$ is an optimal irrigation plan for the measure $\mu$, we conclude
\bel{dIe} \liminf_{\ve\to 0+} {\I^\alpha(\mu_\ve)- \I^\alpha(\mu)\over \ve}
 ~\leq~\lim_{\ve\to 0+} {\E^\alpha(\chi_\ve)- \E^\alpha(\chi)\over \ve}
~=~\int \alpha Z\, d\nu.\eeq
Combining (\ref{dHe}) with (\ref{dIe}) we obtain a contradiction to the optimality
of the solution $(u^*,\mu)$.
\endproof

\begin{remark}{\rm The above argument would fail in higher space dimensions
because, when $d\geq 3$,  a measure $\nu$ supported on a 1-dimensional arc is not absolutely
continuous w.r.t.~capacity.
}\end{remark}

\v

\section{The support of the optimal measure}
\label{s:4}
\setcounter{equation}{0}
Let $(u^*,\mu)$ be an optimal solution to the problem {\bf (OPR)}, describing optimal 
shapes for tree roots.
According to the Theorem~\ref{t:21}, the support of the optimal measure $\mu$ is
contained in the set where the two functions $(1-\psi)u^*$ and $c\alpha Z$ coincide.
In the remainder of this paper, by showing that  these
functions have very different regularity properties, we will prove that the coincidence set
is indeed very small, at least in the case of dimension $d=2$.  We recall that the support of a
positive measure $\mu$ is defined as
$$\Supp(\mu)~\doteq~\Big\{ x\in \R^d\,;~\mu\bigl(B(x, r)\bigr)>0\quad \forall  r>0\Big\}.$$

\begin{theorem}\label{t:41} Let the assumptions {\bf (A1)--(A3)} hold, 
and assume $d=2$, $0<\alpha<1$. 
Let  $(u^*,\mu)$
be an optimal solution to {\bf (OPR)}. 
Then the support of the measure $\mu$ is nowhere dense.
\end{theorem}

Toward a proof, we begin with a few remarks.
\begi
\item[(i)]  The two functions $u^*,\psi$ in (\ref{harveq})-(\ref{adj}) are non-negative and
bounded above.   In particular, there exists constants $K, K'>0$ such that
\bel{Kdef} 0\,\leq\,f(u^*(x))\,\leq \,K,\qquad  \bigl|f'(u^*(x)) \psi(x)\bigr|\,\leq\, K'\qquad\forall x\in \ov \Omega\,.\eeq

\item[(ii)]
By (\ref{harveq}) we have
 \bel{usub}
 \Delta u^*~= u^*\mu - f(u^*)~\geq~-K.\eeq
 Moreover, by (\ref{NC3}) and (\ref{Tus})   it follows
\bel{muoo}\mu\Big(\bigl\{ x\in \ov\Omega\,;~~1-\psi(x)<0\bigr\}\Big)~=~0.\eeq
 Hence  the measure $(1-\psi)\mu$ is non-negative.
By (\ref{adj}) we thus have
 \bel{psup}\Delta\psi~=~-(1-\psi)\mu - f'(u^*) \psi(x)~\leq~K'.\eeq
 
 As a consequence of (\ref{usub})-(\ref{psup}), the function $u^*+K|x|^2$ is
sub-harmonic, while
 $\psi-K' |x|^2$ is super-harmonic.   
In particular (see \cite{AG} for details),  $u^*$ is upper semicontinuous and $\psi$ is lower semicontinuous.
Both $u^*$ and $\psi$ are Borel measurable. Their values are well defined at every 
point $x\in\ov \Omega$. 

\item[(iii)] In addition, since the measure $\mu$ is absolutely continuous w.r.t.~capacity,
we have 
the regularity estimates
\bel{nc9}
u^*\,\in\, H^1(\Omega),\qquad \psi\,\in\, H^1(\Omega),\qquad 
\Phi~\doteq~(1-\psi)\,u^*\,\in\, H^1(\Omega)\cap \L^\infty(\Omega)\,.\eeq
Indeed, from (\ref{Tus}), (\ref{psib}), and the fact that $\psi\geq 0$, we obtain
\bel{Phib}- (\lambda\, u_{max})\cdot u_{max}~\leq~ \Phi(x)~\leq~u^*(x)~\leq~u_{max}\eeq
for all $ x\in\Omega$.
\endi

\v
Combining the previous estimates, we further study the regularity of the product function
$\Phi$ in (\ref{nc9}).
By (\ref{harveq})-(\ref{adj}) it follows
\bel{Dzu}\bega{rl}\Delta\bigl((1-\psi)u^*\bigr)
&=~\ds - u^* \,\Delta \psi + (1-\psi)\,\Delta u^* - 2\nabla \psi\cdot\nabla u^*\\[4mm]
&=~\ds \bigl[ f'(u^*)\psi+(1-\psi)\mu\bigr]u^*  - \bigl[f(u^*)-u^*\mu
\bigr] (1-\psi)- 2\,\nabla \psi\cdot\nabla u^*\\[4mm]
&=~f'(u^*) u^*  \psi - f(u^*)(1-\psi) + 2 (1-\psi) u^*\,\mu  - 2\,\nabla \psi\cdot\nabla u^*.
\enda
\eeq
Therefore, the function $\Phi\doteq (1-\psi)u^*$ provides a bounded solution to the linear elliptic equation
with measure-valued coefficients:
\bel{EM}
\Delta\Phi~=~2\Phi\,\mu+\phi,\eeq
where 
\bel{phi}\phi~\doteq~\bigl[ f'(u^*) u^* + f(u^*)\bigr] \psi - f(u^*)   - 2\,\nabla \psi\cdot\nabla u^*.\eeq
Notice that the product $\Phi\,\mu=(1-\psi)u^*\mu$ is always a positive measure, because of (\ref{muoo}).  However, $\Phi$ can attain both positive and negative values.
We also observe that
$\phi\in \L^1(\Omega)$, because $\psi, u^*\in H^1(\Omega)$.
In addition, at boundary points $x\in \partial\Omega$, the Neumann boundary 
condition holds:
\bel{NC5}\nabla\Phi\cdot\bfn~=~\nabla\bigl((1-\psi)u^*\bigr)\cdot\bfn~=~u^*\,\bigl(\nabla (1-\psi)\cdot \bfn \bigr)+ (1-\psi)\,\bigl(\nabla u^*\cdot \bfn\bigr)~=~0.\eeq
\v
The proof of Theorem~\ref{t:41} will rely on two complementary lemmas.

Given an optimal pair $(u^*,\mu)$,  let $\chi:\Theta\times
\R_+\mapsto \R^2$ be an optimal irrigation plan
for the measure $\mu$. Moreover, consider any particle trajectory
\bel{pat}s\mapsto \gamma(s)\doteq \chi(\bar\theta, s), \eeq for some $\bar\theta\in \Theta$.
By (\ref{stime}), $\gamma$ is  
1-Lipschitz and hence 
a.e.~differentiable.  We denote by
$\bft(s) = \dot\gamma(s)$ the tangent vector.  
We also assume that, on some initial interval, the multiplicity
remains uniformly positive:
\bel{mulg}m(s)~\doteq~|\gamma(s)|\chi~\geq~\delta~>~0\qquad\forall s\in [0,\ov s].\eeq
The first lemma establishes a lower H\"older estimate on the landscape function $Z$.

\begin{lemma}\label{l:42} In the above setting,
for a.e.~$s_0\in [0,\ov s]$
there exist constants $r_0, c_0>0$ such that the following holds.
Calling $x_0=\gamma(s_0)$, the landscape function $Z$ satisfies
\bel{hola}
Z(x)- Z(x_0)~\geq~c_0 |x-x_0|^\alpha\quad\hbox{whenever}\quad 
|x-x_0|< r_0, ~~~\left|\left\langle \bft(s_0)\,,\, {x-x_0\over |x-x_0|}\right\rangle\right|~\leq~{1\over 
3}\,.\eeq
\end{lemma} 

Next, we claim that a converse inequality holds for the function
$\Phi= (1-\psi) u^*$.

\begin{lemma}\label{l:43} In the same setting as Lemma~\ref{l:42}, let 
$\beta\doteq (1+\alpha)/2$, Then for a.e.~$s_0\in [0,\ov s]$
there exists a constant $c_1>0$ and an infinite 
sequence of points $x_k\to x_0 = \gamma(s_0)$ such that 
\bel{holb}
\Phi(x_k)- \Phi(x_0)~\leq~c_1 |x_k-x_0|^\beta,
\qquad \quad \left|\left\langle \bft(s_0)\,,\, {x_k-x_0\over |x_k-x_0|}\right\rangle\right|~\leq~{1\over 
3}\,,\qquad\forall k\geq 1.\eeq
\end{lemma} 

A proof of Lemma~\ref{l:42} will be given in Section~\ref{s:5}, while Lemma~\ref{l:43}
will be proved in Section~\ref{s:6}.
\v
Relying on the two above lemmas, we can now give a proof of Theorem~\ref{t:41}.
Let $y\in \R^2$ be a point inside the support of the measure $\mu$. 
Then, for any given $\ve>0$, there exists  a particle $\bar\theta\in \Theta$
and a path (\ref{pat})  which satisfies (\ref{mulg})  for some $\delta>0$
and such that
$$\bigl|\gamma(\ov s) - y\bigr|~<~\ve.$$

By Lemmas~\ref{l:42} and \ref{l:43}, for a.e.~$s_0\in [0,\ov s]$, at the point $x_0= \gamma(s_0)$  both  (\ref{hola}) and (\ref{holb}) are satisfied.
In particular, we can choose $s_0$ close enough to $\ov s$ so that 
$$|x_0 - y|~\leq~\bigl| \gamma(s_0) - \gamma(\ov s)\bigr|+
\bigl|\gamma(\ov s) - y\bigr|~<~\ve.$$
By (\ref{hola}) and (\ref{holb}), since $\alpha<\beta<1$, there exists a point $x_k$ sufficiently 
close to $x_0$ such that
$$\Phi(x_k)- \Phi(x_0)~\leq~c_1 |x_k-x_0|^\beta~<~c\alpha\cdot c_0 \, |x_k- x_0|^\alpha,
\qquad\qquad |x_k - y|~<~\ve.$$

We now observe that, restricted to the set
$$\Omega^+~\doteq~\bigl\{ x\in\ov\Omega\,;~~1-\psi(x)\,>\,0\bigr\}~=
~\bigl\{ x\in\ov\Omega\,;~~\Phi(x)\,>\,0\bigr\},
$$ the function $\Phi= (1-\psi) u^*$ is the product of two positive, upper semicontinuous 
functions.    Therefore it is upper semicontinuous.   
We can thus find an open neighborhood $V_k$ of $x_k$ such that
\bel{ph5}\Phi(x)- \Phi(x_0)~<~c\alpha\cdot c_0 \, |x- x_0|^\alpha\qquad\forall x\in V_k\,.
\eeq

By Theorem~\ref{t:22} it follows that $\Phi(x_0)\leq c\alpha Z(x_0)$. Together with 
(\ref{ph5})  and (\ref{hola}), this yields
\bel{ph6} \Phi(x)~<~c\alpha\,Z(x_0)+c\alpha\cdot c_0 \, |x- x_0|^\alpha~\leq~
c\alpha Z(x)\qquad\forall x\in V_k\,.\eeq
Hence by Theorem~\ref{t:21}, the open set $V_k$ does not intersect the support of $\mu$.

Since $\ve>0$ was arbitrary,
we conclude that every point $y\in\Supp(\mu)$
lies in the closure of an open set $V$ which does not intersect $\Supp(\mu)$.
This shows that the closed set $\Supp(\mu)$ has empty interior, completing the proof.
\endproof

\v
\section{A lower H\"older estimate for the landscape function}
\label{s:5}
\setcounter{equation}{0}
Aim of this section is to give a proof of Lemma~\ref{l:42}.  Actually, the result
remains valid more generally for  any
positive, bounded Radon measure  $\mu$ on $\R^d$.   For a given $0<\alpha<1$,
let $\chi$ be an optimal irrigation plan for a measure $\mu$, and let $Z$ be the corresponding landscape function.

Let $s\mapsto \gamma(s)$  be a particle trajectory with uniformly positive multiplicity, as in
(\ref{pat})-(\ref{mulg}).
For a given constant $\kappa>0$, consider the set 
\bel{Jk} J_\kappa~\doteq~\Big\{ s_0\in [0,\ov s]\,;~
\bigl| m(s)- m(s_0)\bigr|~\leq~\kappa\, |s-s_0|\quad\forall s\in [0,\ov s]\Big\}.\eeq
Since the multiplicity $m$ is bounded and nonincreasing, 
by Riesz' sunrise lemma
(see for example \cite{KF}, p.~319) it follows
\bel{sunrise}
\meas(J_\kappa) ~\geq~\ov s - {2 m(0)\over \kappa}\,.
\eeq
Therefore 
\bel{sun2}\lim_{\kappa\to +\infty} \meas(J_\kappa) ~=~\ov s.\eeq

In the following, we say that a point $x_0=\gamma(s_0)$ with  $0<s_0<\ov s$ is a {\it good point},
and write $x\in \G$, provided that
\begi
\item[(i)] $s_0$ is a Lebesgue point of the map $s\mapsto \bft(s)$,
\item[(ii)]    $s_0\in J_\kappa$ for some $\kappa$ large enough.
\endi
By (\ref{sun2}) and the fact that $\gamma$ is 1-Lipschitz, 
it follows that the set of good points has full measure.  Namely, 
$\gamma(s_0)\in \G$ for a.e.~$s_0\in [0,\ov s]$.

We claim  that, at every good point $x_0$, the property (\ref{hola}) holds.

\begin{figure}[ht]
\centerline{\hbox{\includegraphics[width=6cm]{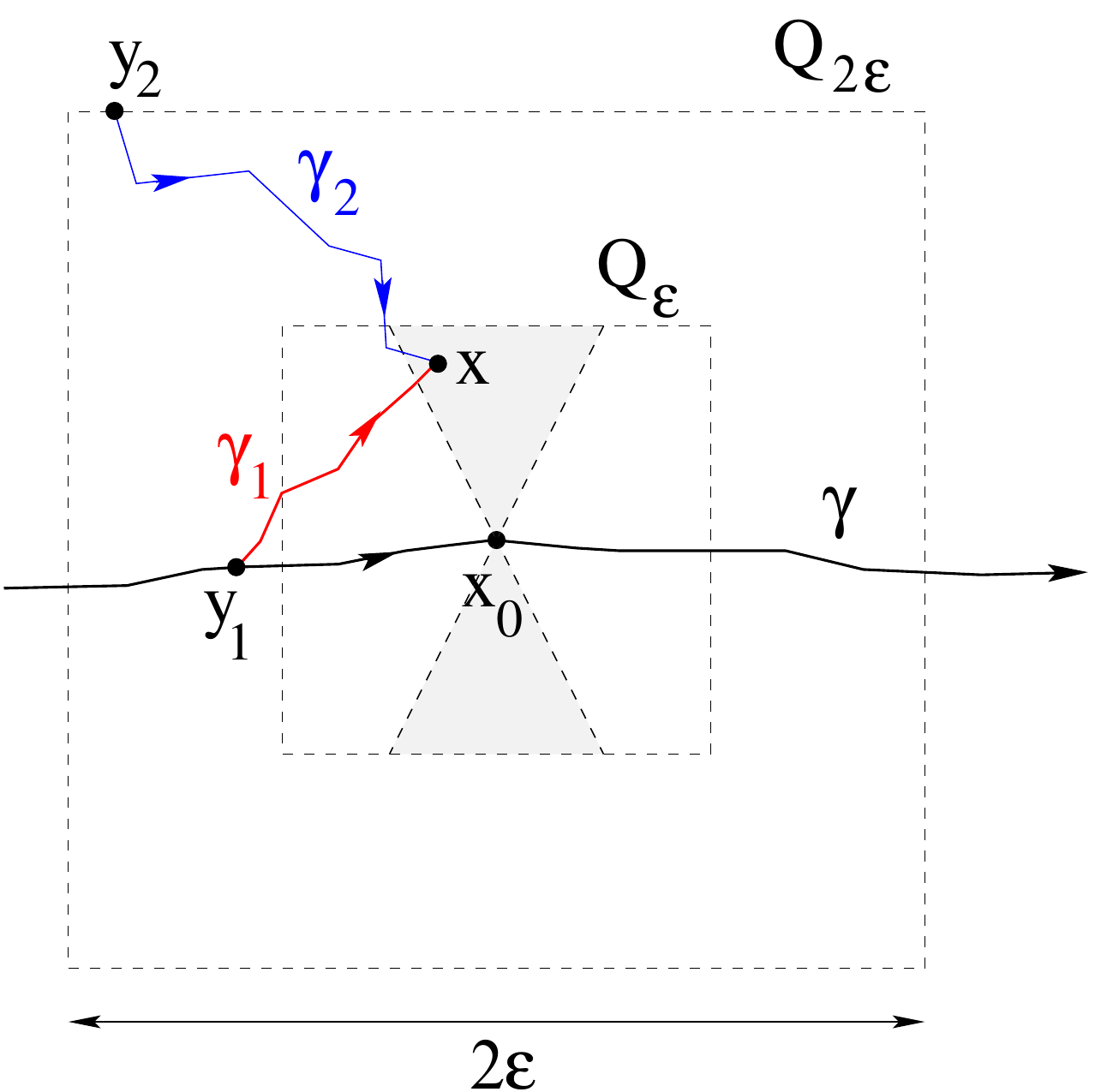}}}
\caption{\small As shown in the proof of Lemma~\ref{l:42}, in the shaded
region inside $Q_\ve$ the landscape function grows at least at a H\"older rate. }\label{f:ir154}
\end{figure}

As shown in Fig.~\ref{f:ir154},
consider a square $Q_\ve$ whose side has length $\ve>0$, centered at $x_0$, with two sides parallel to $\bft(s_0)$.  Let $Q_{2\ve}$ be the concentric square with side of length $2\ve$, and let $N$ be a constant large enough so that 
\bel{Nbig}{N^{1-\alpha}\over 2} - {\sqrt 2\over\alpha}~\geq~1.\eeq
By choosing $\ve>0$ small enough, we can assume that $\gamma$ is the only 
path of multiplicity $> \delta_0/ N$ that intersects $Q_{2\ve}$.

Let $x$ be a point such that
\bel{xqe}x\in Q_\ve\,,\qquad\qquad   \left|\left\langle \bft(s_0)\,,\, {x-x_0\over |x-x_0|}\right\rangle\right|~\leq~{1\over 
3}\,.\eeq
To establish a 
lower bound on $Z(x)$, two cases will be considered.
\v
CASE 1: 
Assume that $x$ is reached by a ``short branch"  $\gamma_1$, that bifurcates
from $\gamma$ at some point
$y_1= \gamma(s_1)\in Q_{2\ve}$.

In this case, at every point along this short branch $\gamma_1$
the multiplicity $m_1$ is bounded by the downward jump in the multiplicity along $\gamma$.
Since $s_0\in \J_\kappa$, this implies
$$m_1~\leq~m(s_1-) - m(s_1+) ~\leq~\kappa \cdot |s_0-s_1 |.$$  
Therefore
\bel{Zx1}Z(x)~\geq~Z(y_1) + \Big(\kappa |s_0-s_1|\Big)^{\alpha-1} |x-y_1|
.\eeq
Since along $\gamma$ the multiplicity is bounded below by (\ref{mulg}), 
we have
\bel{Zy1}
Z(y_1)~\geq~Z(x_0) - \delta_0^{\alpha-1}|s_0-s_1|
.\eeq
The assumption that $\gamma$ is differentiable at $x_0=\gamma(s_0)$
implies that, by choosing $\ve>0$ small enough 
in view of (\ref{xqe}) we can assume
\bel{ss01}|s_0-s_1|~\leq~2|x_0-y_1|,
\qquad\qquad 
\left|\left\langle {x_0-y_1\over |x_0-y_1|}\,,\, {x-x_0\over |x-x_0|}\right\rangle\right|~\leq~{1\over 
2}\,.\eeq
In particular, the angle between these two vectors is larger than $\pi/3$.
By elementary trigonometry, this implies
\bel{xxy}|x-y_1|~\geq~{|x-x_0| \over 2} +{ |x_0-y_1|\over 2}\,.\eeq
Calling 
$$\sigma = |x-x_0|,\qquad r = |x_0-y_1|,$$
from (\ref{Zx1})--(\ref{ss01}) we now obtain
\bel{Z5} \bega{rl} 
Z(x)-Z(x_0)&=~\bigl[ Z(x)-Z(y_1)\bigr] + [Z(y_1) - Z(x_0)\bigr]\\[3mm]
&\geq~\Big(2\kappa |x_0-y_1|\Big)^{\alpha-1} |x-y_1|
-2 \delta_0^{\alpha-1}|x_0-y_1|
\\[3mm]
&\geq~\ds (2\kappa)^{\alpha-1} r^{\alpha-1} {\sigma+r\over 2} - 2\delta_0^{\alpha-1} r
\\[3mm]
&\geq~\ds  2^{\alpha-2}\kappa^{\alpha-1} r^{\alpha-1} \sigma
+ 2^{\alpha-3}\kappa^{\alpha-1} r^\alpha.\enda
\eeq
Notice that the last inequality follows from the fact that $r\leq \sqrt 2\, \ve$, with $\ve>0$ small.
Since the minimum of the right hand side of (\ref{Z5}) is achieved when
$$r~=~{2(1-\alpha)\over \alpha} \sigma,$$
we conclude
\bel{Z6}
Z(x)-Z(x_0)~\geq~c_0\, \sigma^\alpha,
\eeq
for a suitable constant $c_0$.
\v
CASE 2:
Assume that $x$ is reached by a ``long branch"  $\gamma_2$, that enters $Q_{2\ve}$ at a point
$y_2$ which does not lie on the curve $\gamma$.
By (\ref{0xa0}) we have
\bel{Z10}Z(x_0) - Z(y_2)~\leq~{1\over \alpha} \delta_0^{\alpha-1} |x_0-y_2|.\eeq
By construction, this long branch has length $\geq \ve/2$.  Moreover, 
 inside $Q_{2\ve}$, all of  its points have
multiplicity $\leq \delta_0/N$.
In this case, we would have
\bel{Z11}
Z(x)~\geq~Z(y_2) +\left({\delta_0\over N}\right)^{\alpha-1} {\ve\over 2} .\eeq
Together,  (\ref{Z10}) and (\ref{Z11}) yield
\bel{Z12}\bega{rl}
Z(x)- Z(x_0)&\ds \geq~- {1\over \alpha} \delta_0^{\alpha-1} \sqrt 2 \, \ve
+\left({\delta_0\over N}\right)^{\alpha-1} {\ve\over 2}\\[4mm]
&\geq~\ds  \delta_0^{\alpha-1}\ve \left[ {N^{1-\alpha}\over 2} - {\sqrt 2\over\alpha}\right]~\geq~ \delta_0^{\alpha-1}\ve\,, \enda\eeq
because of our choice of  the constant $N$ at (\ref{Nbig}).
Notice that the right hand side of (\ref{Z12}) remains uniformly positive for all $x\in Q_\ve$.

Combining (\ref{Z6}) with (\ref{Z12}), by choosing $r_0<\ve/2$ small enough, we achieve
(\ref{hola}). This achieves the proof of  Lemma~\ref{l:42}.
%
%
%
%
%
%
%
%
%
%
%

\section{Proof of Lemma~\ref{l:43}}
\label{s:6}
\setcounter{equation}{0}
The proof will be worked out in several steps.

{\bf 1.} Let $0<\beta<1$ be given.
To prove the inequality in (\ref{holb}), we need to establish some upper bound 
on the function 
$\Phi$, relying on the fact that it provides a solution to  (\ref{EM}), with $\phi\in \L^1$.   
Since $\mu$ is a positive measure, this can be achieved by constructing a supersolution to 
\bel{EM2} \Delta\Phi~=~-|\phi|,\eeq
with suitable boundary conditions. 
We recall that, by Theorem~\ref{t:22}, along the curve $\gamma$ the upper bound (\ref{49}) holds.   Moreover, at every point $x\in \ov\Omega$, by (\ref{lin}) we have
\bel{EM3}\Phi(x)~\leq~u^*(x)~\leq~u_{max}\,.\eeq
\v
{\bf 2.} As a preliminary, for any integer $\kappa \geq 1$, consider the sets
\bel{Skap}S_\kappa~=~S'_\kappa\cap S''_\kappa\,,\eeq
where
$$S'_\kappa~\doteq~\Big\{ s_0\in [0,\ov s]\,;\quad 
Z(\gamma(s)) - Z(\gamma(s_0))\leq \kappa |s-s_0|\quad\forall s\in [0, \ov s]\Big\},$$
$$S''_\kappa~\doteq~\left\{ s_0\in [0,\ov s]\,;\quad
\limsup_{r\to 0+} {1\over r} \int_{B(\gamma(s_0), r)} \bigl|\phi(x)\bigr|\, dx\,\leq\, \kappa\right\}.
$$
We claim that the union of the sets $S_\kappa$ has full measure in $[0,\ov s]$.
Indeed, since the map $s\mapsto Z(\gamma(s))$ is monotone increasing,
by Riesz' sunrise lemma (see for example \cite{KF}) we have
$$\meas\Big([0,\ov s]\setminus S'_\kappa\Big)~\leq~2   \,{Z\bigl(\gamma(\ov  s)\bigr)-Z\bigl(\gamma(0)\bigr)\over \kappa}\,.$$
This already shows that the union of the sets $S_\kappa'$ has full measure in $[0, \ov s]$.
Next, consider the set
$$\Hat S~\doteq~[0,\ov s]\setminus \bigcup_{\kappa\geq 1} S''_\kappa.$$
If
$ \meas(\Hat S)>0$, 
a contradiction is obtained as follows.
For a given $\kappa>0$, and every $s\in \Hat S$, we can find a sequence of radii
$r_i\downarrow 0$ such that
$$ \int_{B(\gamma(s), r_i)} \bigl|\phi(x)\bigr|\, dx~\geq~\kappa r_i\,.$$
Let $C$ be the constant in (\ref{distlb}).
As $s$ varies in $\Hat S$,
the corresponding intervals $[s- C r_i, s+Cr_i]$ trivially cover $\Hat S$.   
By Vitali's covering theorem (see \cite{EG}), we can extract
a countable family of disjoint intervals $I_j= [s_j- C r_j, s_j+C r_j]$, $j\in \J$, so that the 
collection of intervals $[s_j- 5Cr_j, s_j+5Cr_j]$ covers $\Hat S$.
In particular, this implies
$$\sum_j r_j~\geq~{1\over 10 C} \meas(\Hat S).$$
By (\ref{distlb}), the balls $B(\gamma(s_j), r_j)$ are mutually disjoint.
Hence
\bel{70}\sum_{j\in\J} \int_{B(\gamma(s_j), r_j)}   |\phi|\, dx  ~\leq~\|\phi\|_{\L^1}\,.\eeq
On the other hand,
\bel{71}\sum_{j\in\J} \int_{B(\gamma(s_j), r_j)}   |\phi|\, dx  ~\geq~\sum_j \kappa r_j~\geq~ 
\kappa\cdot {1\over 10 C} \meas(\Hat S).\eeq
Since $\kappa$ can be arbitrarily large, if $\meas(\Hat S)>0$, from 
(\ref{71})  we obtain a contradiction
with (\ref{70}).

 We now define a point $x= \gamma(s)$ to be a {\em good point} if the
tangent vector $\bft(s) = \dot \gamma(s)$ is well defined, and
if $s\in S_\kappa$ for some $\kappa\geq 1$.
In the remainder of the proof, we will show that (\ref{holb}) holds for every good point $x_0$.
\v
{\bf 3.} Toward a future comparison, we study two elliptic problems
on a domain $\D_\delta\subset\R^2$ which, in polar coordinates, has the form
\bel{Ddelta}\D_\delta
~\doteq~\bigg\{(r,\theta)\,;~~r\in [0,1],~~\theta\in \Big[-{\pi\over 2}-\delta\,,~{\pi\over 2}+
\delta\Big]\bigg\},\eeq
for some $\delta > 0 $ small.
On this domain, let
$\Phi_1$ be the solution to
\bel{Ph1}\left\{ \bega{rll} \Delta\Phi_1(x)&=~0 \qquad &\hbox{if}\quad x\in \D_\delta\,,\\[3mm]
\Phi_1(x)&\ds=~ |x| \qquad &\hbox{if}\quad x\in \partial \D_\delta\,.
\enda\right.\eeq
In addition, we consider the solution $\Phi_2$  to the Poisson problem
\bel{Ph2} \left\{ \bega{rll}
\Delta\Phi_2 &=~-\phi \qquad &\hbox{on}~~ \D_\delta,\\[3mm]
\Phi_2 &=~0\qquad &\hbox{on}~~ \partial \D_\delta\,,\enda\right.\eeq
assuming that $\phi\geq 0$ and, for some $\kappa \geq 1$,
\bel{ipbo}
\int_{\D_\delta\cap B(0,r)}\phi\, dx~\leq~\kappa r\qquad\qquad\forall r\in [0,1].\eeq
Both problems (\ref{Ph1})-(\ref{Ph2}) are more conveniently studied by constructing a conformal map, transforming the domain   $\D_\delta$ into the half disc
$$\D'~=~\bigg\{(r,\theta)\,;~~r\in [0, 1],~~
\theta\in \Big[-{\pi\over 2}\,,~{\pi\over 2}\Big]\bigg\}.$$

\begin{figure}[ht]
\centerline{\hbox{\includegraphics[width=11cm]{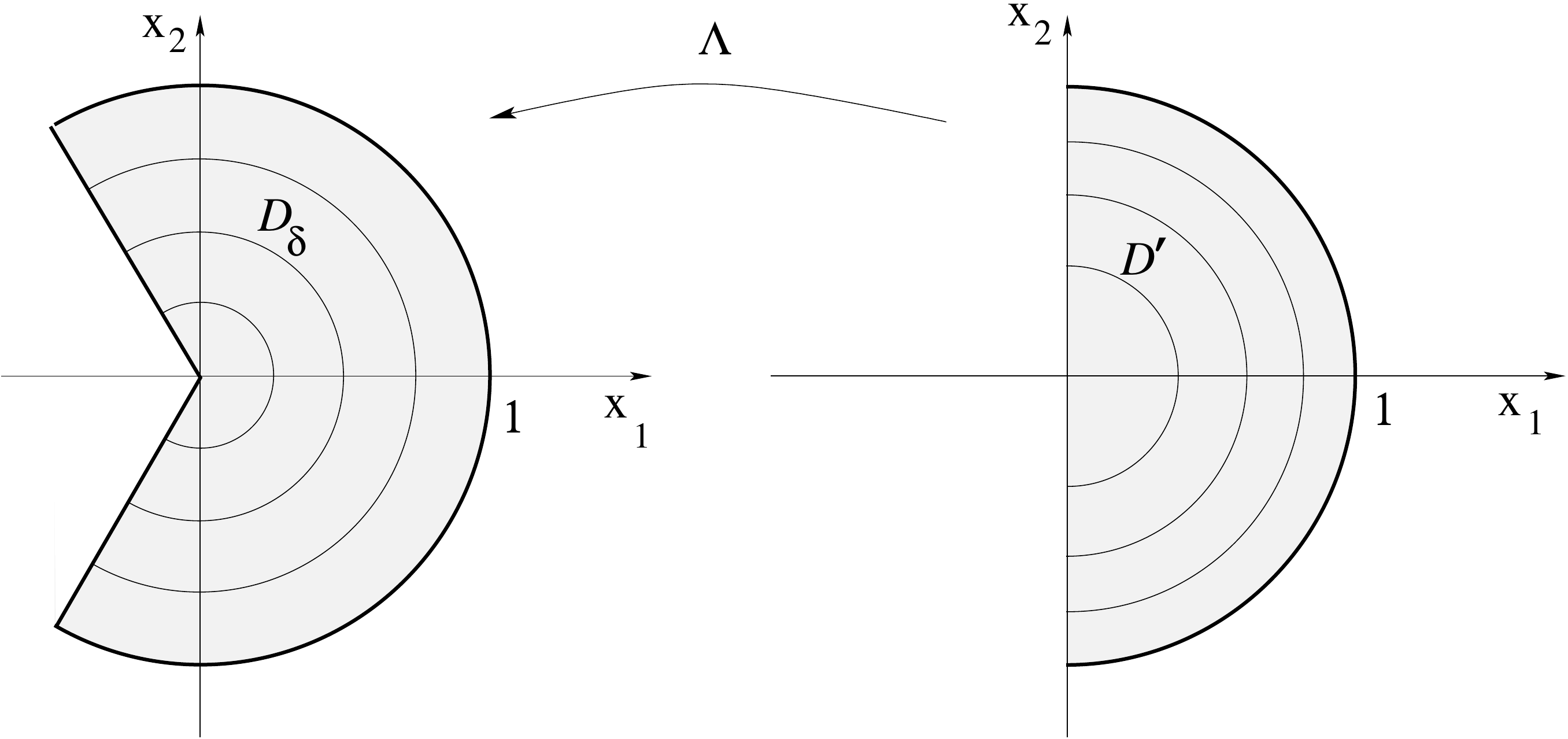}}}
\caption{\small  Taking  $p = 1+{2\delta\over\pi}$, the conformal map $z\mapsto z^p$
transforms  the half disc $D'$ into the domain  $D_\delta$.
}\label{f:ir148}
\end{figure}

 Choosing $p= 1 + {2\delta \over \pi}$, 
the transformation 
$z\mapsto\Lambda(z)=z^p$ in the complex plane is then a conformal map
from $D'$ onto  $\D_\delta$,
as  shown in Fig.~\ref{f:ir148}.
In polar coordinates, this takes the form
\bel{conft}
(\tilde r, \tilde\theta)~=~\Lambda(r,\theta)~=~(r^p, p \theta).\eeq
\v
{\bf 4.} To construct a solution to (\ref{Ph1}),
 we now solve the problem on the half disc
\bel{DP2}
\left\{ \bega{rll} \Delta v&=~0\qquad\qquad& x\in \D',\\[3mm]
v&=\ds ~r^p\qquad\qquad & x\in \partial \D'.\enda \right.\eeq
The function
\bel{Udef}
\Phi_1(z)~\doteq~v(z^{1/p})\eeq will then provide the desired solution to (\ref{Ph1}).

In the following, we shall not need an explicit expression for $\Phi_1$, but only 
an estimate on its asymptotic behavior near the origin.
To construct $v$, we can use the Green's function on the half space,
then add a correction (smooth in a neighborhood of the origin) to take into account the
effect of the boundary at $|z|=1$.  Restricted to the subdomain
$$\D''~\doteq~\left\{ x=(x_1,x_2)\,;~|x|<{1\over 2}\,,~ x_1>0\right\},$$
this leads to
\bel{vrep} v(x_1,x_2) ~= ~ { x_1 \over \pi} \int_{-1}^1 { |y_2|^p \over x_1^2 + (y_2-x_2)^2 } dy_2
+ e(x_1, x_2),\eeq
for some smooth correction term $e$, with $e(0,0)=0$.

To estimate $v(x)$ we shall rely on the following identity, valid for $x_1 > 0$.
  \bel{DPint}
    \frac{x_1}{\pi}\int_{-\infty}^{\infty} \frac{1}{x_1^2 + (t-x_2)^2}\,dt~ = ~1.
  \eeq
 The integral in (\ref{vrep}) can be bounded as 
\bel{I12}\bega{l} \ds
 \frac{x_1}{\pi}\int_{-1}^1  \frac{|y_2|^p}{x_1^2 + (y_2-x_2)^2}\,dy_2 \\[4mm]
 \qquad\ds = \frac{x_1}{\pi}  \int_{-2|x|}^{2|x|}  \frac{|y_2|^p}{x_1^2 + (y_2-x_2)^2}\,dy_2 + \frac{x_1}{\pi}  \int_{\{2|x| < |y_2| \le 1 \}} \frac{|y_2|^p}{x_1^2 + (y_2-x_2)^2}\,dy_2\\[4mm]
 \qquad \doteq~I_1 + I_2\,. \enda\eeq
By (\ref{DPint}) we now have
\bel{I1} I_1~\leq~2^p |x|^p.\eeq
To estimate the second integral, we observe that $|y_2| > 2|x|$ implies
  $$ |y_2| ~\le~ \bigl|(0,y_2) - (x_1,x_2)\bigr| + |x| < \sqrt{x_1^2 + (y_2-x_2)^2} + \frac12 |y_2|. $$
  Therefore $x_1^2 + (y_2-x_2)^2 > \frac14 |y_2|^2$, and hence
  $$ \frac{x_1}{\pi}  \int_{\{2|x| < |y_2| \le 1 \}} \frac{|y_2|^p}{x_1^2 + (y_2-x_2)^2}\,dy_2 ~<~ \frac{4 x_1}{\pi}  \int_{\{2|x| < |y_2| \le 1 \}} |y_2|^{p-2} \,dy_2.$$
We thus conclude that, if $|x|<1/2$, the second integral can be estimated as
 \bel{I2}I_2~\leq~ \frac{4 x_1}{\pi}  \int_{\{2|x| < |y_2| \le 1 \}} |y_2|^{p-2} \,dy_2 ~=~ \frac{8 x_1}{\pi}  \int_{2|x|}^1 t^{p-2} \,dt~ =~ 
    \frac{8 x_1}{(p-1)\pi} \bigl(1-(2|x|)^{p-1}\bigr).\eeq
Furthermore,   since the correction term $e$ is smooth, it can be bounded above by some linear function:  $e(r,\theta) \le C r$.
 
 Here and throughout the following, for notational convenience we denote by $C>0$
a positive constant, whose value can change at each step.

\begin{figure}[ht]
\centerline{\hbox{\includegraphics[width=10cm]{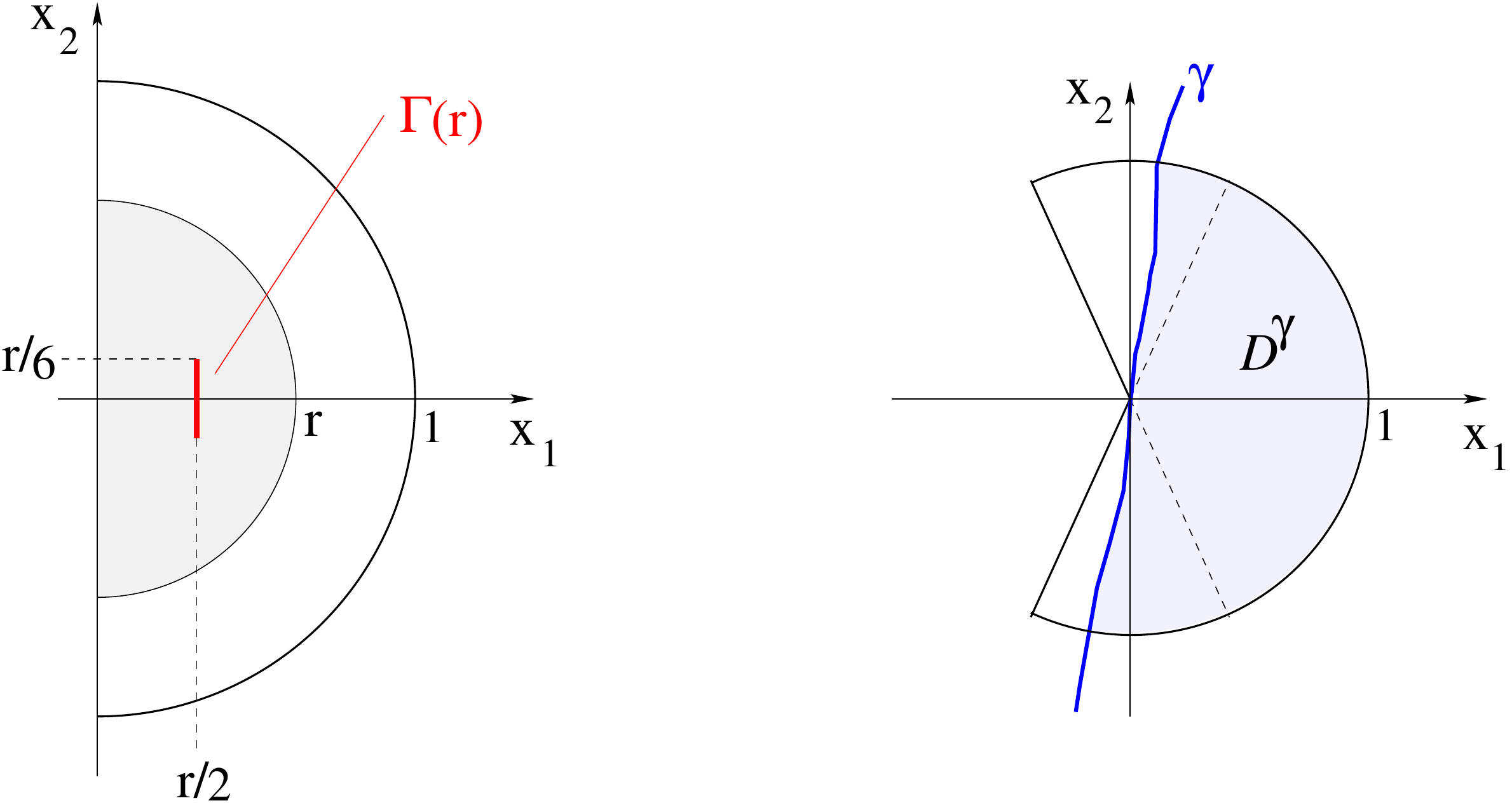}}}
\caption{\small Left: the interval $\Gamma(r)$, where the average
value  for the solution $u$ of (\ref{PPu}) can be estimated. 
Right: the domain $\D^\gamma$ considered at (\ref{Ph12}).}\label{f:ir155}
\end{figure}

Combining (\ref{I1})-(\ref{I2}), we obtain the estimate
\bel{vlin}
\bigl|v(x)\bigr|~\leq~C\, |x|\qquad\qquad\forall  x\in \D',\eeq
for a suitable constant $C$.  In turn, by (\ref{Udef}), this implies
\bel{phi1}
\bigl|\Phi_1(x)\bigr|~\leq~C |x|^{1/p}.\eeq
\v
{\bf 5.} In addition to the upper bound (\ref{phi1}),
we observe that the solution $\Phi_1$ of (\ref{Ph1}) satisfies the lower bound
\bel{P1up}\Phi_1(x)~\geq~|x|\qquad \qquad x\in \D_\delta\,.\eeq
Indeed, one immediately checks that the function $\vp(x)=|x|$ is a subsolution
to (\ref{Ph1}).
\v
{\bf 6.} We now consider the  solution to (\ref{Ph2}).
Using polar coordinates, if  $\Phi_2=\Phi_2(\tilde r, \tilde\theta)$ is a solution to (\ref{Ph2})
on $\D_\delta$, then
 the function 
$u(r,\theta)~=~\Phi_2(r^p, p\theta)$
satisfies 
$$\bega{rl}\Delta u(r,\theta)
&\ds =~u_{rr} + {1\over r} u_r + {1\over r^2} u_{\theta\theta}\\[3mm]
&=~p^2 r^{2p-2}\Phi_{2,rr} + p^2 r^{p-2} \Phi_{2,r}
+p^2 r^{-2} \Phi_{2,\theta\theta}\\[3mm]
&\ds=~p^2 r^{2p-2} \Delta \Phi_2( r^p , p \theta).
\enda
$$
We thus set 
$$f(r,\theta)~\doteq~p^2 r^{2p-2} \phi(r^p, p\theta),$$
and consider the Poisson problem on the half disc
\bel{PPu}
\left\{ \bega{rll}\Delta u &=~-f\qquad &\hbox{on}~\D',\\[3mm]
u&=~0\qquad  & \hbox{on}~\partial \D'.\enda\right.\eeq
Notice that
$$\|f\|_{\L^1(\D')}~=~\int_{-\pi/2}^{\pi/2} \int_0^1
p^2 r^{2p-1} \phi(r^p, p\theta) dr d\theta~=~\int_{\D_\delta} 
p^2 r^{2p-1} \phi(\tilde r,\tilde\theta) {d\tilde r\over p r^{p-1}} \, {d\tilde\theta
\over p}  ~=~\|\phi\|_{\L^1(\D_\delta)}\,.
$$
We observe that, since the function $f$ is only in $\L^1$,  pointwise bounds on $u$
cannot be deduced from a Sobolev embedding theorem.  However, we can establish 
a bound on the average value of $u$ on an interval $\Gamma(r)$, as shown in Fig.~\ref{f:ir155}, left.

Fix a radius $0<r\leq 1$ and consider the  interval
\bel{Gr}\Gamma(r)~\doteq
~\left\{(x_1, x_2)\,;\quad x_1  = {r\over 2},\quad |x_2| < {r\over 6}  \right\},\eeq
shown in Fig.~\ref{f:ir155}, left.
Setting 
$$y=(y_1,y_2),\qquad y'=(-y_1, y_2),$$
an upper bound for the solution $u$ of (\ref{PPu}) will be obtained by using the Green's formula
for the half space $\R^2_+ \doteq \{(x_1,x_2)\,;~ x_1 \geq 0\}$.   
For a given radius $r>0$, 
it will be convenient to split the function $f$ as
$$f(x) = f^\flat + f^\sharp~=~f\cdot {\bf 1}_{\{ |x|\leq r\}} +f\cdot {\bf 1}_{\{ |x|> r\}}\,.
$$
This leads to
\bel{usf}\bega{rl}u(x)&\leq~u^+(x)\ds\doteq
~{1\over 2\pi} \int_{\R^2_+} \bigl( \ln |x-y'| - \ln|x-y|\bigr) f(y)\, dy\\[4mm]
&\ds=~{1\over 2\pi} \int_{\R^2_+} \bigl( \ln |x-y'| - \ln|x-y|\bigr) \bigl(f^\flat(y)+ f^\sharp(y)\bigr)\, dy
~\doteq~u^\flat(x) + u^\sharp(x)\,.
\enda\eeq
In the next steps, we shall prove an integral bound on $u^\flat$ and a pointwise bound on $u^\sharp$.
\v
{\bf 7.}
Setting
$$z_2\,=\, x_2-y_2\,,\qquad F^\flat (y_1)~=~\int_{\R} f^\flat(y_1,y_2)\, dy_2\,,$$
\bel{ax0} \Tilde G(x_1, y_1)~\doteq~{1\over 4\pi} \int_{\R} \Big( \ln\bigl((x_1 + y_1)^2+z_2^2 \bigr) - 
\ln\bigl((x_1 - y_1)^2+ z_2^2\bigr)\Big)\, dz_2\,,\eeq
we now compute
\bel{ax1}\bega{rl}U^\flat (x_1)
&\ds \doteq ~~{1\over 2\pi} \int_{\R} \int_{\R^2_+}  \bigl( \ln |x-y'| - \ln|x-y|\bigr) f^\flat (y)\, dy\,dx_2
\\[4mm]
&=\ds~{1\over 4\pi} { \int_{\R}\int_{\R}\int_0^\infty } \Big( \ln\bigl((x_1 + y_1)^2+ z_2^2\bigr) - 
\ln\bigl((x_1 - y_1)^2+ z_2^2\bigr)\Big) f^\flat (y_1,y_2) dy_1 dy_2\, dz_2\\[4mm]
&=\ds~ {\int_0^\infty} \Tilde G(x_1, y_1)\, F^\flat (y_1) dy_1\,.
\enda
\eeq
 From the representation (\ref{ax0}), since $x_1\geq 0$, we have 
	\bel{ax00} \Tilde G(x_1, y_1)~\geq~0 \, ,\quad \Tilde G(x_1, - y_1)~=~ - \Tilde G(x_1, y_1)\, , \qquad\forall y_1~\geq~ 0\, . \eeq
We now claim that, when $r >0$ is sufficiently small, one has
\bel{ax2} U^\flat (x_1)~\leq~C x_1 |\ln x_1|\cdot  \|f^\flat\|_{\L^1} . \eeq
Indeed, by (\ref{ax1}) we can write 
\bel{ax3}\bega{rl}\ds U^\flat(x_1)&\ds=~\int_{0\leq y_1\leq 2x_1} \Tilde G(x_1, y_1) F^\flat(y_1)\, dy_1 + \int_{ y_1> 2x_1} \Tilde G(x_1, y_1) F^\flat(y_1)\, dy_1\\[3mm]
 &\ds=~\int_{0\leq y_1\leq 2x_1} \Tilde G(x_1, y_1) F^\flat(y_1)\, dy_1,
 \enda
   \eeq
because $2x_1 = r$ and $F^\flat(y_1)$ vanishes for $|y_1| > r$.
By the definition (\ref{ax0}) it follows
   \bel{ax4}\bega{rl}\ds \Tilde G(x_1, y_1)&\ds=~{1\over 4\pi} \int_{\R} \Big( \ln(z_2^2 +(x_1 + y_1)^2) - 
   \ln(z_2^2 +(x_1 - y_1)^2)\Big)\, dz_2 \\[3mm]
   &\ds=~{1\over 2\pi} \int_{0}^{|y_1|} \Big( \ln(z_2^2 +(x_1 + |y_1|)^2) - 
   \ln(z_2^2 +(x_1 - |y_1|)^2)  \Big)\, dz_2 \\[3mm]
   &\ds \quad +~ {1\over 2\pi} \int_{|y_1|}^\infty \ln \left( 1+ {4x_1|y_1|\over z_2^2 + (x_1 - |y_1|)^2} \right)\, dz_2\\[3mm]
   &\ds \doteq~A + B.
     \enda\eeq
A direct computation yields 
     \bel{ax5} |A|~\leq~{1\over 2\pi} \int_0^{|y_1|} |\ln(z_2^2 + (x_1 + |y_1|)^2)|\, dz_2
      + {1\over 2\pi} \int_0^{|y_1|} |\ln(z_2^2 + (x_1 - |y_1|)^2) |\, dz_2\, .  \eeq
     Choosing $r$ small enough, since we have $|y_1|\leq 2x_1 = r$, in (\ref{ax5}) we can assume
     $$ z_2^2 + (x_1 + |y_1|)^2 ~\leq~ 1\, ,\quad  z_2^2 + (x_1 - |y_1|)^2~\leq~1.$$
  This yields the estimate 
     \bel{ax6} |A|~\leq~{1\over \pi} \int_0^{|y_1|} |\ln(z_2^2)|\, dz_2~\leq~{2\over \pi} 
     \bigl(|y_1|\ln(|y_1|) + |y_1| \bigr)~\leq~C x_1|\ln(x_1)|\, . \eeq
     To estimate $B$ in (\ref{ax4}), using the inequality  $\ln (1 + s)\leq s $ for $s\geq 0$, we
     obtain
     \bel{ax7} |B|~\leq~ {1\over 2\pi} \int_{|y_1|}^\infty {4x_1|y_1|\over z_2^2}\, dz_2~=~{4x_1\over 2\pi}~\leq~C x_1|\ln(x_1)|\, . \eeq
     {}Using (\ref{ax4}) together with the bounds (\ref{ax6}) and (\ref{ax7}),
     from  (\ref{ax3}) we conclude  that, for $r > 0$ sufficiently small,
	 \bel{ax8} U^\flat(x_1)~\leq~\int_{0\leq y_1 \leq 2x_1} (|A| + |B|) 
	 F^\flat(y_1)\, dy_1~\leq~ C x_1 |\ln (x_1)| \int_{0\leq y_1 \leq 2x_1} F^\flat(y_1)\,  dy_1\, . \eeq
\v
\v{\bf 8.} It now remains to estimate
\bel{ush}
u^\sharp(x)~\doteq~{1\over 2\pi} \int_{\R^2_+} \bigl( \ln |x-y'| - \ln|x-y|\bigr) f^\sharp(y)\, dy\eeq
for $x$ on the segment $\Gamma(r)$ at (\ref{Gr}).  
Setting 
$$F^\sharp(s)~\doteq~\int_{r<|y|<s} f^\sharp(y)\, dy~\leq~\kappa s,$$
we compute
\bel{Fss} \bega{rl}\ds
\int_{r<|y|<1} { f^\sharp (y)\over |y|}\, dy&\ds=~\int_r^1 {1\over s} \cdot \left({d\over ds} F^\sharp(s)\right)
\, ds~=~\left[{ F^\sharp(s)\over s}\right]_r^1 + \int_r^1 {1\over s^2} \cdot F^\sharp(s)
\, ds\\[4mm]
&\ds\leq~\kappa + \int_r^1 {\kappa \over s}\, ds ~\leq~\kappa\bigl(1+|\ln r|\bigr).
\enda
\eeq
Next, we observe that $x \in \Gamma(r)$ implies  $|x| \leq {\sqrt{10} \over 6}r$. This  leads to the bound
$$ \ln |x-y'| - \ln|x-y| ~=~ {1\over2}\ln\biggl( 1+ {4x_1 y_1\over |x-y|^2} \biggr) \le \frac{2x_1 y_1}{(|y|-|x|)^2} ~< ~{2 |x| |y| \over (|y|/3)^2}. $$
Using this bound in (\ref{ush}),  for $x\in \Gamma(r)$  we obtain
\bel{ush2}
u^\sharp(x)~\leq~C\,  |x| \int_{|y|>r} { f^\sharp (y)\over |y|}\, dy~\leq~
C r\, \kappa\bigl(1+|\ln r|\bigr),\eeq
for some constant $C$.
\v
{\bf 9.}  Since $u^+ = u^\flat + u^\sharp$,
combining the integral bound (\ref{ax2}) with the pointwise bound (\ref{ush2})
we obtain
\bel{uub}\int_{\Gamma(r)} u^+(x)\, dx~=~
\int_{-r/6}^{r/6} u^+\left({r\over 2}, \, x_2\right)\, dx_2~\leq~C r |\ln r|\cdot  \|f^\flat\|_{\L^1}  +{r\over 2} \cdot
 C r  \kappa\bigl(1+|\ln r|\bigr).\eeq
We now recall that, by  the assumption (\ref{ipbo}),
$$ \|f^\flat\|_{\L^1}~=~\int_{|x|<r} f(x)\, dx~=~\int_{|x|< r^p} \phi(x)\, dx ~\leq~C r^p,$$
for some constant $C$ and $p>1$.
Using this inequality in (\ref{uub}) and  (\ref{usf}), we conclude that  
the average value of $u$ over $\Gamma(r)$ satisfies the  bound 
\bel{avva}
\avint_{\Gamma(r)} u\, dx~\leq~\avint_{\Gamma(r)} u^+\, dx~\leq~C r\bigl(1+|\ln r|\bigr),\eeq
for a suitable constant $C$ and all $r>0$ sufficiently small.
\v
{\bf 10.}  
Next, consider the more general problem
\bel{Ph12} \left\{ \bega{rl}
\Delta\Phi (x)&=~-\phi(x), \qquad x\in \D^\gamma,\\[3mm]
\Phi (x)&=~|x|,~\qquad x\in \partial \D^\gamma\,.\enda\right.\eeq
As shown in Fig.~\ref{f:ir155}, right, the domain $\D^\gamma$ is the portion of the unit disc
to the right of a Lipschitz curve $\gamma$, with  
\bel{HSd}\gamma~\subset~ \Hat \D_\delta~\doteq~\bigg\{ (r,\theta)\,;~r\in [0,1],~~\theta\in 
\left[ -{\pi\over 2} -\delta\,,~-{\pi\over 2} +\delta\right]\cup\left[ {\pi\over 2} -\delta\,,~{\pi\over 2} +\delta\right]\bigg\}.\eeq
As before, we assume that  (\ref{ipbo}) holds, for some $\kappa \geq 1$. 
Moreover, since we are seeking an upper bound on the solution $\Phi$ of (\ref{Ph12}), 
w.l.o.g.~we can assume that $\phi\geq 0$.   
To cover the general case it suffices to replace $\phi$ by its positive part
$\phi^+(x)\doteq \max\{ \phi(x), 0\}$.
 
In view of (\ref{P1up}) and the fact that $\Phi_2 \geq 0$, we have the comparison
\bel{PPP}
\Phi(x)~\leq~\Phi_1(x)+\Phi_2(x)\qquad\qquad  \forall x\in \D^\gamma,\eeq
where $\Phi_1$ and $\Phi_2$ are the solutions to (\ref{Ph1}) and (\ref{Ph2}), respectively.
\v
{\bf 11.}
Thanks to (\ref{avva}), we can now construct a sequence of points
$P_k\in \Gamma(r_k)$, with $r_k\to 0$, 
such that 
$$u(P_k)~\leq~u^+(P_k)~\leq~C  |P_k|\Big(1+\bigl|\ln |P_k|\bigr|\Big).$$
Setting $Q_k = \Lambda (P_k) = (P_k)^p$ we now obtain
\bel{p2b}\Phi_2(Q_k)~\doteq~u(P_k)~\leq~C  |P_k|\Big(1+\bigl|\ln |P_k|\bigr|\Big)~
\leq~C |Q_k|^{1/p}\Big(1+\bigl|\ln |Q_k|\bigr|\Big).\eeq
We recall that $C$ always denotes a positive constant, whose precise value may change
at each occurrence.

Given $0<\beta<1$, we can now choose $\delta>0$ small enough so that 
$${1\over p} ~=~\left[ 1 + {2\delta\over \pi}\right]^{-1} ~>~\beta.$$
Combining (\ref{phi1}) with (\ref{p2b}),  we thus obtain
\bel{PQk}\Phi(Q_k)~\leq~\Phi_1(Q_k) + \Phi_2(Q_k)~\leq~C |Q_k|^{1/p}+ C |Q_k|^{1/p}\Big(1+\bigl|\ln |Q_k|\bigr|\Big)~<~C  |Q_k|^\beta.\eeq
By (\ref{Gr}),  the assumption $P_k=(P_{k1}, P_{k2})\in \Gamma(r_k)$ implies
$$|P_{k2}|~\leq~{P_{k1}\over 3}\,.$$
Hence, calling $\bfe_2\doteq (0,1)$, we have
$$ \left|\left\langle \bfe_2\,,\, {P_k\over |P_k|}\right\rangle\right|~=~
{|P_{k2}|\over\sqrt{P_{k1}^2 + P_{k2}^2}} ~\leq~{1\over \sqrt{10}}~<~{1\over 
3}\,.$$
In turn, if $\delta>0$ was chosen sufficiently small, then $Q_k=(Q_{k1}, Q_{k2})=\Lambda(P_k)$ still satisfies
\bel{Qk3} \left|\left\langle \bfe_2\,,\, {Q_k\over |Q_k|}\right\rangle\right|~\leq~{1\over 
3}\,.\eeq
\v
{\bf 12.} At last, we can now complete the proof of Lemma~\ref{l:43}.
Let $x_0=\gamma(x_0)$ be a {\it good point}, as defined at the end of step {\bf 2.}
By a possible rotation of coordinates, we can assume that the tangent vector is
$\bft(s_0) = \bfe_2 = (0,1)$.
We then choose a small radius $\rho>0$ and consider an affine
transformation 
$$x~\mapsto~ y ~=~ \T x$$
mapping the
disc $B(x_0, \rho)$ centered at $x_0$ with radius $\rho$ onto the disc $B(0,1)$
centered at the origin with unit radius. 

Restricted to the disc $B(x_0,\rho)$, the function $\Phi$ satisfies the elliptic equation
(\ref{EM}) together with the lower bound
$$\Phi(x)~\leq~c\alpha Z(x)\qquad\qquad \hbox{for}~~x\in\gamma.$$
Since $x_0$ is a {\it good point},  we can choose $\lambda, \rho>0$ small enough 
so that
the corresponding function
$$\Tilde \Phi(y)~=~\lambda\bigl[\Phi(\T^{-1} y) - \Phi(x_0)
\bigr]\qquad\qquad \hbox{for}~~|y|\leq 1$$
satisfies  a system of the form
\bel{Phhh} \left\{ \bega{rll}
\Delta\Tilde \Phi (y)&\geq~-\phi(y), \qquad & y\in \D^\gamma,\\[3mm]
\Tilde \Phi (y)&\leq~|y|,~\qquad & y\in \partial \D^\gamma\,.\enda\right.\eeq
In other words, $\Tilde\Phi$ provides a subsolution to 
(\ref{Ph12}).  
By the previous analysis, there exists an infinite sequence of points $Q_k\to 0$
such that 
$$\Tilde \Phi(Q_k)~\leq~C\, |Q_k|^\beta.$$
Going back to the original coordinate $x\in B(x_0,\rho)$, this yields a sequence of points
$q_k = \T^{-1} Q_k$ such that 
$$q_k\to x_0,\qquad\qquad \Phi(q_k)-\Phi(x_0)~\leq~C\, |q_k-x_0|^\beta.$$
Notice that the inequality 
\bel{ie1} \left|\left\langle \bft(s_0)\,,\, {q_k-x_0\over |q_k-x_0|}\right\rangle\right|~\leq~{1\over 
3}\eeq
is an immediate consequence of (\ref{Qk3}).
This completes the 
 proof.
\endproof
\v
{\bf Acknowledgments.}
The research of A.~Bressan and Q.~Sun
was partially supported by NSF  with grant DMS-1714237, ``Models of controlled biological growth".
The research of S.~T.~Galtung was partially supported by the grant ``Wave Phenomena and Stability -- a Shocking Combination (WaPheS)'' (project no.\ 286822) from the Research Council of Norway.
 \v

\end{document}